\documentstyle{amltd2004}
\begin{document}
\input amssym.def
\input amssym.tex

\def\lrp#1{\left(#1\right)}
\annalsline{156}{2002}
\received{June 30, 2000}
\startingpage{213}
\def\bye{\end{document}}
 \font\tenrm=cmr10

\catcode`\@=11
\font\twelvemsb=msbm10 scaled 1100
\font\tenmsb=msbm10
\font\ninemsb=msbm10 scaled 800
\newfam\msbfam
\textfont\msbfam=\twelvemsb  \scriptfont\msbfam=\ninemsb
  \scriptscriptfont\msbfam=\ninemsb
\def\msb@{\hexnumber@\msbfam}
\def\Bbb{\relax\ifmmode\let\next\Bbb@\else
 \def\next{\errmessage{Use \string\Bbb\space only in math
mode}}\fi\next}
\def\Bbb@#1{{\Bbb@@{#1}}}
\def\Bbb@@#1{\fam\msbfam#1}
\catcode`\@=12

 \catcode`\@=11
\font\twelveeuf=eufm10 scaled 1100
\font\teneuf=eufm10
\font\nineeuf=eufm7 scaled 1100
\newfam\euffam
\textfont\euffam=\twelveeuf  \scriptfont\euffam=\teneuf
  \scriptscriptfont\euffam=\nineeuf
\def\euf@{\hexnumber@\euffam}
\def\frak{\relax\ifmmode\let\next\frak@\else
 \def\next{\errmessage{Use \string\frak\space only in math
mode}}\fi\next}
\def\frak@#1{{\frak@@{#1}}}
\def\frak@@#1{\fam\euffam#1}
\catcode`\@=12

\newcommand{\lap}{\ensuremath{\triangle}} 


 \title{Superlogarithmic estimates on\\ pseudoconvex domains and CR manifolds}

 \shorttitle{Superlogarithmic estimates} 

 \acknowledgements{Research was partially supported by NSF Grant
DMS-9801626.}
 \author{J. J. Kohn}
 \institutions{Princeton University,
 Princeton, NJ\\
{\eightpoint {\it E-mail address\/}: kohn@princeton.edu
}}

\centerline{\bf Abstract}
\vglue9pt

This paper is concerned with proving superlogarithmic estimates for the operator
$\Box_b$ on pseudoconvex CR manifolds and using them to establish hypoellipticity of
$\Box_b$ and of the $\bar{\partial}$-Neumann problem. These estimates are established
under the assumption that subellipticity degenerates in certain specified ways.

\section{Introduction}

 Let $\Omega\subset X$ be a domain with a smooth boundary in a complex manifold~$X$. If $h$ is a holomorphic function
on $\Omega$ which is smooth up to the boundary we denote by $\dot h$ the restriction of $h$ to the boundary of
$\Omega$. Then $\dot h$ satisfies the the tangential Cauchy-Riemann equations: $\bar\partial_b\dot h=0$. The operator
$\bar\partial_b$ can be extended to forms on the boundary; this extension is useful in studying the boundary values of
$\bar\partial$-closed forms (see [KR]). In [K1] I defined $\bar\partial_b$ and the associated Laplacian
\begin{equation}
\Box_b=\bar\partial_b\bar\partial_b^*+\bar\partial_b^*\bar\partial_b 
\end{equation}
on forms on abstract CR manifolds. The operator
$\Box_b$ is then used to study $\bar\partial_b$. The analysis of $\Box_b$ on $(p,q)$-forms on compact strongly pseudoconvex CR
manifolds $M$ of dimension greater than three with $0<q<1/2({\rm dim} M -1)$ is carried out in [K1]. To analyze the operator
$\Box_b$ on $(p,q)$-forms with $0\le q\le 1/2 ({\rm dim} M -1)$ and $ {\rm dim} M\ge 3$ we need to use microlocal
analysis; see [K3]. Microlocal analysis is also used to prove that the operators $\bar\partial_b$ and $\Box_b$ have a
closed range if $M$ is the boundary of a domain in a Stein manifold (see [K4]). In this paper we show how microlocal
estimates for $\Box_b$  imply estimates for the $\bar\partial$-Neumann problem. In particular we investigate
superlogarithmic estimates but the  technique presented here works for subelliptic and other types of estimates as well.

Let $M$ be a pseudoconvex ${\rm CR}$ manifold of dimension $2n-1$.  We say that $\Box_b$ is
hypoelliptic at $x_o\in M$ if whenever $\Box_b\varphi$ is $C^\infty$ in some neighborhood of $x_o$ then
$\varphi $ is $C^\infty$ in a neighborhood of $x_o$. If $\Box_b$ is hypoelliptic on $(0,1)$ forms then we
conclude that: 1. The restriction of $\bar{\partial}_b$ to the orthogonal complement of 
square integrable CR functions (i.e. functions in the null space of $\bar{\partial}_b$) is
hypoelliptic. 2. If the square integrable function $f$ is $C^\infty$ in a neighborhood     
of $x_o\in M$ then $S_bf$, the orthogonal projection of $f$ to the space of square-integrable CR functions, is also $C^\infty$
in the same neighborhood. 

Hypoellipticity  is a consequence of subellipticity (see [KN]). The operator $\Box_b$ is
\emph{subelliptic} at $x_o$ if there exist a neighborhood $U$ and positive constants  $\varepsilon$ and
$C$ such that
\begin{equation}
\| \varphi\|^2_\varepsilon \le C(\Box_b\varphi,\varphi) = CQ_b(\varphi,\varphi),
\end{equation}
for all $\varphi\in C_0^\infty (U)$.
Such estimates for the closely related $\bar{\partial}$-Neumann problem have been studied
extensively. The necessary and sufficient condition for such an estimate to hold for the
$\bar{\partial}$-Neumann problem near a point in the boundary of a pseudoconvex domain is that the D'Angelo type at the
point is finite (see [D] and [C]). The approach to the problem taken in [K2] is via subelliptic multipliers; this leads to the
condition of ``finite ideal type", which is sufficient for subellipticity.    This condition is probably equivalent to
the  finite D'Angelo type. In the appendix there is a discussion of these matters and of their extension to CR manifolds.

Throughout this paper we will deal with $(p,q)$-forms. However, in all proofs we will restrict
ourselves to the case of $(0,q)$-forms, since we consider only local properties where the extension
from type $(0,q)$ to type $(p,q)$ is trivial. To prove {\it a~priori} estimates for $(0,q)$-forms we
will write the decomposition $\varphi =\varphi^0+\varphi^++\varphi^-$; the precise definition will
be given in the next section. The method is to prove estimates analogous to (1.2) with $\varphi$
replaced by $\varphi^0$, $\varphi^+$, and $\varphi^-$, respectively. The case of $\varphi^0$ being
particularly simple, we have
\begin{equation}
\| \varphi^0\|^2_1 \le CQ_b(\varphi^0,\varphi^0);
\end{equation}
this is the ``elliptic'' case with $\varepsilon =1$ and holds on all CR manifolds. To study the other
microlocalizations we consider the ideals ${\cal I}^+_q(x_0)$ and  ${\cal I}^-_q(x_0)$ defined
(approximately) as follows (the precise definition will be given in Section 5). The
ideal ${\cal I}^+_q(x_0)$ consists of all germs of $C^\infty$ functions at $x_0$ such that
$\rho\in{\cal I}^+_q(x_0)$ if
\begin{equation}
\| \rho\varphi^+\|^2_\varepsilon \le
C\bigl(Q_b(\varphi^+,\varphi^+)+\|\varphi^+\|^2),
\end{equation}
and $\rho\in{\cal I}^-_q(x_0)$ if
\begin{equation}
\|\rho \varphi^-\|^2_\varepsilon \le
C\bigl(Q_b(\varphi^-,\varphi^-)+\|\varphi^-\|^2).
\end{equation}
Thus if $1\in {\cal I}^+_q(x_0)\cap{\cal I}^-_q(x_0)$ we obtain (1.2). In the appendix we
give explicit formulas for elements of these ideals and also the definition of finite ideal
$q$-type. In general we have ${\cal I}^+_q(x_0)\subset {\cal I}^+_{q'}(x_0)$ if $q\le q'$ and
${\cal I}^+_q(x_0)={\cal I}^-_{n-q-1}(x_0)$. Hence,
$$
{\cal I}^+_q(x_0)\cap{\cal I}^-_q(x_0)=\cases{{\cal I}^+_q(x_0) & if
$q\le\frac{1}{2}(n-1)$\cr{\cal I}^-_q(x_0) & if $q\ge\frac{1}{2}(n-1)$}.
$$  

We also have ${\cal I}^+_0(x_0)={\cal I}^-_{n-1}(x_0)=\{0\}$ and so the estimate (1.2) cannot
hold for $(0,0)$-forms, i.e. functions, or for $(0,n-1)$-forms. In fact, for pseudoconvex compact CR
manifolds, the nullspaces of $\Box_b$ on functions and on $(0,n-1)$-forms are in general infinite-dimensional (on
boundaries of domains in complex manifolds  $\Box_bh=0$ is equivalent to the condition that $h$ is the restriction of a
holomorphic function). It is then natural to impose the global condition of orthogonality to the nullspace of $\Box_b$ but
even this assumption is not sufficient. Rossi (see [R]) has given an example of a strongly pseudoconvex three-dimensional
CR manifold for which the Szeg\"o projection does not preserve smoothness and hence the restriction
of $\Box_b$ to the orthogonal complement of its nullspace is not hypoelliptic. In [K3] it is shown
that for compact pseudoconvex CR manifolds on which the range of $\bar{\partial}_b$ is closed in
$L_2$ and for which $1\in{\cal I}^+_1(x_0)$ the restriction of $\Box_b$ to the orthogonal
complement of its nullspace is globally hypoelliptic. This closed-range property is proved for compact CR
manifolds that are boundaries of domains in Stein manifolds (see [K4]).

The problem is to find conditions under which hypoellipticity holds when subellipticity fails. The
most desirable condition is a purely geometric condition such as in Fedii's theorem (see [F]), which
states that the operator $\frac{\partial^2}{\partial x^2} + a(x)\frac{\partial^2}{\partial y^2}$ with
$a\in  C^\infty$ is hypoelliptic if $a(x)>0$ when $x\ne 0$.  In [K6] Fedii's result is generalized to
degenerate subelliptic operators in any dimension. This phenomenon appears on CR manifolds;  thus
on a pseudoconvex CR manifold $\Box_b$ is hypoelliptic in a neighborhood $U$ of $x_o$ if it is
subelliptic at all $x$ with $x\ne x_o$. The natural geometric condition for hypoellipticity seems to
be that the  points where subellipticity fails are contained in a real curve transversal to the
``good'' directions (i.e. the tangent vectors to the curve do not lie in $T^{(1,0)}+T^{(0,1)}$); this
is true for certain classes of CR manifolds (see [K5] and [K6]) but Christ has shown that it is false in general
(see [Ch2]).  The example of Kusuoka and Stroock (see [KS]) illustrates that, in general, the geometry
alone of the set on which the operator degenerates does not give sufficient information to determine
whether it is hypoelliptic. They study the operator $\frac{\partial^2}{\partial x^2} +
a^2(x)\frac{\partial^2}{\partial y^2} +\frac{\partial^2}{\partial z^2}$ in ${\Bbb R}^3$ with
$a(x)\ne 0$ if $x\ne 0$ and prove that it is hypoelliptic if and only if $\lim_{x\to 0}x\log a(x)=0$. 
Hypoellipticity of such operators is obtained by use of superlogarithmic estimates (see [Ch1] and [M]).
Superlogarithmic estimates for $\Box_b$ are formulated as follows.
 
\numbereddemo{Definition}
If $M$ is a CR manifold and $x_o\in M$ then a {\it superlogarithmic estimate} holds 
for $\Box_b$ on positively microlocalized $(p,q)$-forms at $x_o$ if there exists a neighborhood $U$
of $x_o$ such that for each $\delta >0$ there exists $C_\delta$ such that
\begin{equation}
\|(\log\Lambda)\varphi^+\|^2\le \delta Q_b(\varphi^+ ,\varphi^+ ) + C_\delta\|\varphi^+\|^2 ,
\end{equation}
for all $(p,q)$ forms $\varphi $ with support in $U$. Here $\Lambda $ denotes the square root of the
Laplacian. A {\it superlogarithmic estimate} holds for $\Box_b$ on negatively
microlocalized $(p,q)$-forms at $x_o$ if the above holds with $\varphi^+$ replaced by $\varphi^-$.
In that case 
\begin{equation}
\|(\log\Lambda)\varphi^-\|^2\le \delta Q_b(\varphi^- ,\varphi^- ) + C_\delta\|\varphi^-\|^2.
\end{equation}
\enddemo

Here we will prove superlogarithmic estimates of the form (1.6) and (1.7) when there exists a
subelliptic multiplier $\rho$ which may go to zero on a submanifold $S$ at a controlled rate. The
condition which we will use is the following:
\begin{equation}  
\lim_{x\to S}d(x,S)\log|\rho (x)| = 0,
\end{equation}
where $d(x,S)$ denotes the distance from $x$ to $S$. Our results will impose conditions on the {\it
holomorphic dimension} of $S$; this concept will be recalled in Section 4.

\proclaim{Theorem}
Let $M$ be a pseudoconvex {\rm CR} manifold. Let $S\subset M$ be a manifold such that the holomorphic
dimension of $S$ at each point is less than or equal to $q-1$. Let $x_0\in S$ and $\rho\in
{\cal I}^+_q(x_0)${\rm .} Suppose that $\rho$ satisfies {\rm (1.8);} then there exists a neighborhood $U$ of
$x_0$ on which the superlogarithmic estimate {\rm (1.6)} is satisfied for all $(0,q)$\/{\rm -}\/forms with support
in $U${\rm .}
\endproclaim

\proclaim{Theorem}
Let $M$ be a pseudoconvex {\rm CR} manifold. Let $S\subset M$ be a manifold such that the holomorphic
dimension of $S$ at each point is less than or equal to $n-q-2$. Let $x_0\in S$ and $\rho\in
{\cal I}^-_q(x_0)$. Suppose that $\rho$ satisfies {\rm (1.8);} then there exists a neighborhood
$U$ of $x_0$ on which the superlogarithmic estimate {\rm (1.7)} is satisfied for all $(0,q)$\/{\rm -}\/forms with support
in $U${\rm .}
\endproclaim

One of the crucial steps in proving these results is the following localization lemma, which may be
of independent interest.

 \proclaimtitle{Localization Lemma}
\proclaim{Lemma}
Let $M$ be pseudoconvex {\rm CR} manifold and let $S\subset M$ be a submanifold of holomorphic dimension less
than or equal to $m$. Let $x_0\in S$. For small positive $a$ let $S_a$ denote the set of $x\in M$
such that $d(x,S)<a$. Then there exists a neighborhood $U$ of $x_0$ and a constant $C$ such that 
\begin{equation}
\| \varphi^+ \|^2_{S_a} \le C\Bigl( a^2Q_b(\varphi^+, \varphi^+) + 
\| \varphi^+ \|^2_{M -S_a}+\|\varphi^+\|^2_{-1}\Bigr),
\end{equation}
for all $(p,q)$\/{\rm -}\/forms $\varphi$ with $q\ge m+1$ supported in $U$. Furthermore{\rm ,}
\begin{equation}
\| \varphi^- \|^2_{S_a} \le C\Bigl( a^2Q_b(\varphi^-, \varphi^-) + 
\| \varphi^- \|^2_{M -S_a}+\|\varphi^+\|^2_{-1}\Bigr),
\end{equation}
for all $(p,q)$-forms $\varphi$ with $q\le n-m-2$ supported in $U$. Here $\|\ \|_{S_a}$\ and
$\|\ \|_{M-S_a}$ denote $L_2$-norms over $S_a$ and $M-S_a${\rm ,} respectively.    
\endproclaim

These theorems imply that the superlogarithmic estimate
\begin{equation}   
\|(\log\Lambda)\varphi\|^2\le \delta Q_b(\varphi ,\varphi) + C_\delta\|\varphi\|^2
\end{equation}
holds, when $M$ is pseudoconvex, for $(0,q)$-forms supported in $U$ whenever (1.8) with 
$\rho\in{\cal I}_q^+(x_o)\cap {\cal I}_q^-(x_o)$ holds for $S$ of holomorphic dimension
$m\le\break\min\{q-1,n-q-2\}$. Note that the estimate (1.11) cannot hold if $q=0$ or $q=n-1$.

The superlogarithmic estimates (1.6), (1.7), and (1.11) imply the following results concerning
hypoellipticity.

\proclaim{Theorem}
 Let $M$ be a pseudoconvex {\rm CR} manifold and let $x_0\in M$ such that {\rm (1.11)}
 holds for $(p,q)$\/{\rm -}\/forms supported in a
neighborhood $U$ of $x_0$. Suppose that $\varphi$ and $\alpha$ are $(p,q)${\rm -}forms on $M$ in $L_2(M)$
such that $\Box_b(\varphi)=\alpha$. Suppose further that the restriction of $\alpha$ to $U$ is in
$C^\infty(U)$. Then the restriction of $\varphi$ to $U$ is also in $C^\infty(U)$.
\endproclaim

 Denote by ${\cal H}^{p,q}_b(M)$ the nullspace of $\Box_b$ on $(p,q)$-forms.  If $M$ is
compact and if (1.11) holds for $(p,q)$-forms in some neighborhood of every point of $M$ then ${\cal H}^{p,q}_b(M)$
consists of $C^\infty$ forms and is finite-dimensional. Furthermore, under this assumption, we
have the following hypoellipticity consequences:

\vglue2pt
(a) Whenever $\alpha\in L_2(M)$ and $\alpha\perp{\cal H}^{p,q}_b(M)$ the equation
$\Box_b(\varphi)=\alpha$ has a unique solution $\varphi\perp{\cal H}^{p,q}_b(M)$. 

\vglue2pt
(b) If $\bar{\partial}_b\alpha =0$ then there exists a unique $(p,q-1)$-form $\psi$ such that
$\bar{\partial}_b\psi =\alpha$ and $\psi$ is perpendicular to the nullspace of $\bar{\partial}_b$;
the restriction of $\psi$ to any   open set on which $\alpha$ is $C^\infty$ is also $C^\infty$.

\vglue2pt
(c) Dually, if $\bar{\partial}^*_b\alpha =0$ then there exists a unique $(p,q+1)$-form $\theta$
such that $\bar{\partial}^*_b\theta =\alpha$ and  $\theta$ is perpendicular to the nullspace of
$\bar{\partial}^*_b$; the restriction of $\theta$ to any   open set on which $\alpha$ is
$C^\infty$ is also $C^\infty$.

\vglue2pt
(d) The orthogonal projection operators on ${\cal H}^{p,q}_b(M)$ and on the nullspaces of
$\bar{\partial}_b$ and $\bar{\partial}_b$ are pseudolocal (that is, the the operator applied to a
form is $C^\infty$ on the open set on which the form is $C^\infty$).      

\vglue2pt
The next theorem deals with the cases of $(p,n-1)$ and $(p,0)$-forms. To deal with these we will
assume that the range of $\bar{\partial}_b$ on functions is closed. We will also assume that the
estimate (1.6) \pagebreak holds for $(p,n-1)$-forms which is equivalent to (1.7) for $(p,0)$-forms. In
view of the above results these conditions are satisfied when $M$ is pseudoconvex and there
exists $\rho\in {\cal I}_1^-(x_0)$ which satisfies (1.8) with $S$ of holomorphic dimension
zero (i.e. $S$ totally real).

\proclaim{Theorem}   
Let $M$ be a {\rm CR} manifold such that the range of $\bar{\partial}_b$ on functions is closed. Let
$x_0\in M$ such that {\rm (1.6)} holds for $(p,n-1)$\/{\rm -}\/forms and {\rm (1.7)} holds for $(p,0)$\/{\rm -}\/forms supported 
in a neighborhood $U$ of $x_0$. Suppose that $\varphi$ and $\alpha$ are $(p,n-1)${\rm -}forms{\rm ,} on $M$ in
$L_2(M)${\rm ,} which are orthogonal to ${\cal H}_b^{p,n-1}(M)${\rm ,} such that $\Box_b(\varphi)=\alpha$.
Suppose further that the restriction of $\alpha$ to $U$ is in $C^\infty(U)$. Then the restriction of
$\varphi$ to $U$ is also in $C^\infty(U)$. Dually{\rm ,} the same holds when $\varphi$ and $\alpha$ are
$(p,0)$\/{\rm -}\/forms{\rm ,} on $M$ in $L_2(M)${\rm ,} which are orthogonal to ${\cal H}^{p,0}_b(M)${\rm ,} such that
$\Box_b(\varphi)=\alpha$. 
\endproclaim

Now the spaces ${\cal H}^{p,0}_b(M)$ and ${\cal H}^{p,n-1}_b(M)$ are, in general, infinite
dimensional. Nevertheless, if (1.5) holds for $(p,n-1)$-forms at each point of $M$ we have the 
properties:
 
\vglue3pt
(a) The closed range of $\bar{\partial}_b$ on functions implies that $\Box_b$ has closed range on
all forms. Hence, for $q=0$ and $q=n-1$, whenever $\alpha\in L_2(M)$ and
$\alpha\perp{\cal H}^{p,n-1}_b(M)$ the equation $\Box_b(\varphi)=\alpha$ has a unique solution
$\varphi\perp{\cal H}^{p,n-1}_b(M)$. 

\vglue3pt
(b) If $\alpha$ is a $(p,n-1)$-form   orthogonal to ${\cal H}^{p,n-1}(M)$ with
$\bar{\partial}_b\alpha =0$ then there exists a unique $(p,n-2)$-form $\psi$ such that
$\bar{\partial}_b\psi =\alpha$ and $\psi$ is perpendicular to the nullspace of $\bar{\partial}_b$;
the restriction of $\psi$ to   the open set on which $\alpha$ is $C^\infty$ is also $C^\infty$.

\vglue3pt
(c) Dually, if $\alpha$ is a $(p,0)$-form orthogonal to ${\cal H}_b^{p,0}(M)$
with $\bar{\partial}^*_b\alpha =0$ then there exists a unique $(p,1)$-form $\theta$ such that
$\bar{\partial}^*_b\theta =\alpha$ and  $\theta$ is perpendicular to the nullspace of
$\bar{\partial}^*_b$; the restriction of $\theta$ to   the open set on which $\alpha$ is
$C^\infty$ is also $C^\infty$.

\vglue3pt
(d) The orthogonal projection operators on ${\cal H}_b^{p,0}(M)$, ${\cal H}_b^{p,n-1}(M)$ and on
the nullspaces of $\bar{\partial}_b$ in $(p,n-2)$-forms and $\bar{\partial}_b$ in $(p,1)$-forms 
are pseudolocal.

\vglue3pt

Finally we take up the applications of superlogarithmic estimates to the $\bar{\partial}$-Neumann
problem. Let $\Omega$ be a domain in a complex $n$-dimensional manifold $X$. Assume that the
boundary of $\Omega$, denoted by $b\Omega$, is $C^\infty$. For $(p,q)$-forms on $\Omega$ we
define the operator $\Box$ by
\begin{equation}
\Box = \bar{\partial}\bar{\partial}^*+\bar{\partial}^*\bar{\partial},
\end{equation}
where $\bar{\partial}^*$ is the $L_2$ adjoint of $\bar{\partial}$. Then $\Box\varphi$ is a form
whose components are Laplacians of the corresponding components of $\varphi$ which
satisfies the $\bar{\partial}$-Neumann boundary conditions: $\varphi$ is contained in the domain
of $\bar{\partial}^*$ on $(p,q)$-forms and $\bar{\partial}\varphi$ is contained in the domain
of $\bar{\partial}^*$ on $(p,q+1)$-forms. Let ${\cal H}^{p,q}(\Omega)$ denote the null space of
$\Box$. Then we have:

\proclaim{Theorem}
Suppose that $x_0\in b\Omega$ and suppose that {\rm (1.6)} holds for $(p,q)$-forms on $b\Omega$
supported in $U\cap b\bar{\Omega}${\rm ,} where $U$ is a neighborhood of $x_0$ in $X$ and $1\le q\le
n-2$. Suppose further that $\Box\varphi =\alpha$ with $\alpha$ a $(p,q)$\/{\rm -}\/form in $L_2$ whose
restriction to $U\cap \bar{\Omega}$ is in $C^\infty$. Then the restriction of $\varphi$ to $U\cap
\bar{\Omega}$ is in~$C^\infty$.   
\endproclaim

Note that the hypothesis of the above theorem is satisfied if $b\Omega$ is pseudoconvex and if
there exists $\rho\in{\cal I}^+_q(x_o)$ satisfying (1.8) with $S$ of holomorphic dimension
less than or equal to $q-1$. The theorem has consequences which are   analogous to those of
Theorem  1.5 listed above.

\proclaim{Theorem}
Suppose that $x_0\in b\Omega$ and suppose that {\rm (1.6)} holds for  $(p,n-1)$\/{\rm -}\/forms on $b\Omega$
supported in $U\cap b\bar{\Omega}${\rm ,} where $U$ is a neighborhood of $x_0$ in $X$. Suppose that $X$
is a Stein manifold and that $\Omega$ is pseudoconvex with $\bar{\Omega}$ compact. Suppose
further that $\Box\varphi =\alpha$ with $\alpha$ a $(p,n-1)$\/{\rm -}\/form in $L_2$ whose restriction to
$U\cap\bar{\Omega}$ is in $C^\infty$. Then{\rm ,} if $\varphi\perp{\cal H}^{p,n-1}(\Omega)${\rm ,} the
restriction of
$\varphi$ to $U\cap\bar{\Omega}$ is in $C^\infty$.   
\endproclaim

Note that the hypotheses of this theorem are satisfied if there exists 
$\rho\in{\cal I}^+_{n-1}(x_o)$ satisfying (1.8) with $S$ of holomorphic dimension less than or
equal to $n-1$. Again, the consequences of this theorem are entirely analogous to those of Theorem
1.6 listed above.

Here I wish to express my thanks to M. Christ; communication with him led to the research presented here.
I am also indebted to him for supplying the example discussed at the end
of Section 4. I also thank Siqi Fu and the referee 
for going through the original manuscript and suggesting several improvements and clarifications.

\vglue-12pt
 
\section{Definitions and preliminaries}
\vglue-6pt
 
 \advance\theoremcount by 1
{\it Definition} 2.1.
Let $M$ be a $2n-1$ dimensional differentiable manifold.
Denote by ${\Bbb C}{\rm T}(M)$ the complexified tangent bundle of $M$. An {\it integrable}  
 CR {\it structure} on $M$ is given by a subbundle
$T^{1,0} \subset {\Bbb C} T(M)$ with the following properties. If $T^{0,1}$
denotes the conjugate of $T^{1,0}$ then $T_x^{1,0}\cap T_x^{0,1} = \{ 0 \}$ for all $x$. The 
subbundle $T^{1,0} \oplus T^{0,1}$ is of codimension one. If $L$ and $L'$ are
local vector fields with values in $T^{1,0}$ then the commutator $[L,L'] =LL'-L'L$
also has values in  $T^{1,0}$. We say $M$ is a  CR {\it manifold} if $M$ has
a given integrable CR structure.

\numbereddemo{Definition}
A {\it form of type $(p,q)$} at $x \in M$ is a skew-symmetric multilinear map
$$
\varphi : \underbrace{{T_x^{1,0}} \times \cdots\times {T_x^{1,0}}}_{p}\times
\underbrace{{T_x^{0,1}} \times \cdots \times {T_x^{0,1}}}_{q} \to {\Bbb C}.
$$  
The {\it bundle of $(p,q)$\/{\rm -}\/forms} is denoted by ${\cal A}_b^{p,q}$.
\enddemo 

\numbereddemo{Definition}
The operator $\bar{\partial}_b :{\cal A} ^{p,q}_b \to {\cal A} ^{p,q+1}_b$ is 
defined as follows. If $\varphi \in {\cal A} ^{p,q}$ let $\varphi '$ be a $(p,q)$-form
which restricted to $\prod_p T^{1,0}\times \prod_q T^{0,1}$
equals $\varphi$. Then $\bar{\partial}_b\varphi$ is the restriction of $d\varphi '$ to
$\prod_p T^{1,0}\times \prod_{q+1} T^{0,1}$. The operator $\bar{\partial}^*_b :{\cal A}^{p,q}_b
\to
{\cal A} ^{p,q-1}_b$ is the $L_2$-adjoint of $\bar{\partial}_b$. 
\enddemo 

{\it Definition} 2.4.
Denote by $\theta$ a nonvanishing real $1$-form which annihilates $T^{1,0} \oplus T^{0,1}$.
Then the {\it Levi form} on $M$ is the hermitian form on $T^{1,0}$ given by 
$\sqrt{-1}d\theta (L,\bar{L'})$, where $L$ and $L'$ are in $T^{1,0}$. Now,  $M$ is {\it pseudoconvex}
if for some choice of $\theta$ the Levi form is nonnegative.  
\vglue6pt \advance\theoremcount by 1

Let $L_1,\ldots ,L_{n-1}$ be a local basis for $(1,0)$ vectorfields in a neighborhood $U$
of $x_o \in M$ and let $\omega_1,\ldots ,\omega_{n-1}$ be the dual basis of $(1,0)$-forms.
If $u \in C^\infty (U)$,
$$
\bar{\partial}_b(u) = \sum_j \bar{L}_j(u)\bar{\omega}_j.
$$
\vglue-6pt\noindent
If $\varphi=\sum\varphi_j\bar\omega_j$ is a $(0,1)$-form on $U$,
\begin{eqnarray*}
\noalign{\vskip-6pt}
\bar{\partial}_b\varphi& = &
\sum_{i<j}\Bigl(\bar{L}_j{\varphi}_i - \bar{L}_i{\varphi}_j +
\sum_ka^k_{ij}\varphi_k\Bigr)\bar{\omega}_i
\wedge\bar{\omega}_j
\\  \noalign{\vskip-6pt}
\noalign{\noindent 
and} \noalign{\vskip-3pt}
\bar{\partial}^*_b\varphi& =& -\sum_i(L_i(\varphi_i) +a_i\varphi_i),\\
\noalign{\vskip-16pt}
\end{eqnarray*}
where $a^k_{ij},\  a_i \in C^\infty (U)$.

In general the operator $\bar{\partial}_b :{\cal A} ^{0,q}_b \to {\cal A} ^{0,q+1}_b$ is expressed
as follows. Let $\varphi \in {\cal A} ^{0,q}_b$ be the form given locally by:
\vglue4pt\centerline{${\displaystyle
\varphi = \sum_I \varphi_I \bar{\omega}_I ,
}$}
\vglue4pt\noindent
where $I$ is the $q$-tuple of integers $I=(i_1,\ldots ,i_q)$ with $0<i_1<\cdots <i_q \le n-1$ and
where $\bar{\omega}_I = \bar{\omega}_{i_1} \wedge \cdots\wedge \bar{\omega}_{i_q}$.  Then 
$$
\bar{\partial}_b\varphi = \sum_K\lrp{\sum_{j\notin I}\epsilon^{jI}_K\bar{L}_j\varphi_I+a_{IK}\varphi_I} \bar\omega_K,
$$
where $K$ runs through all strictly increasing $(q+1)$-tuples and each of the coefficients $\epsilon^{jI}_K$ is either $0, 1,$
or
$-1$ and is defined as follows. First, if $j\notin I$ we denote by $\langle jI\rangle $ the ordered
$q$-tuple containing $j$ and the elements of $K$. Then we define
$$
\epsilon^{jI}_K=\cases{0 & if $j\notin K$\cr
\mathop{{\rm sgn}}\nolimits\left({jI\atop K}\right) & if $j\in K$,}
$$
where $\mathop{{\rm sgn}}\nolimits\left({jI\atop K}\right)$ denotes the sign of the permutation
$\{j,I\}\to K$.
\noindent
Further,
$$
\bar{\partial}_b^*\varphi = -\sum_H\lrp{\sum_{i,I\supset H}\epsilon^{iH}_I L_i\varphi_I
+a_H\varphi_H}\bar{\omega}_H,
$$ 
where $H$ runs through all strictly increasing $(q-1)$-tuples and 
$a_{IK},\  a_H\in C^\infty(U)$.

We choose real coordinates $\{x_1,\ldots x_{2n-1}\}$ with the origin at $x_o$ such 
 that, \ \ $\theta (\frac{\partial}{\partial x_{2n-1}})>0$ and such that, setting $z_i = x_i
+\sqrt{-1}x_{i+n-1}$ for $i=1,\ldots ,n-1$,  we have $L_i|_{x_o} = \frac{\partial}{\partial z_{i}}|_{x_o}$.
Set
$T=-\sqrt{-1}\frac{\partial}{\partial x_{2n-1}}$. Then  the Levi form can be written as $c_{ij}
=d\theta (L_i,\bar{L}_j)$ and we have
$$
[L_i,\bar{L}_j] = c_{ij}T +\sum_k(d^k_{ij}L_k +e^k_{ij}\bar{L}_k),
$$
where $d^k_{ij},\  e^k_{ij} \in C^\infty (U)$.
 
Let $\{\xi_1,\ldots ,\xi_{2n-1}\}$ be the dual coordinates to
$\{x_1,\ldots ,x_{2n-1}\} $ and  $|\xi|^2 = \sum_i{\xi^2_i}$.
 Let $\psi^+$ and $\tilde{\psi}^+$ be   nonnegative functions in 
$C^\infty (\{\xi\in{\Bbb R}^{2n-1}\Big{|} |\xi| =1\})$, with range in $[0,1]$, such that 
\begin{eqnarray*}
&&{\rm supp}(\psi^+) \subset \{|\xi|=1\  \big{|}\  \xi_{2n-1}\ge
\frac{1}{2}| \xi '|\},
\\
&&{\rm supp}(\tilde{\psi}^+) \subset \{|\xi|=1\  \big{|}\  \xi_{2n-1}\ge
\frac{1}{4}| \xi '|\},
\end{eqnarray*}
$\psi^+(\xi)=1$ when $\xi_{2n-1}\ge \frac{3}{4}| \xi'|$ and $\tilde\psi^+(\xi)=1$ when 
$\xi_{2n-1}\ge \frac{1}{3}| \xi '|$.  Here  $\xi'=(\xi_1,\ldots \xi_{2n-2})$. We extend these functions to
${\Bbb R}^{2n-1}$ so that $\psi^+ (\xi ) = \psi^+(\frac{\xi}{|\xi|})$ and $\tilde\psi^+ (\xi ) =
\tilde\psi^+(\frac{\xi}{|\xi|})$, for $|\xi | \ge 1$ and so that $\psi^+(\xi)=\tilde\psi^+(\xi)=0$ for
$|\xi |<\frac{1}{2}$ with $\tilde\psi^+(\xi)=1$ on ${\rm supp}(\psi^+)$. Set
$\psi^-(\xi) = \psi^+(-\xi)$, $\tilde\psi^-(\xi) = \tilde\psi^+(-\xi)$ and
$\psi^0 = 1 -\psi^+ -\psi^-$. Define $\tilde\psi^0$ so
that $\tilde\psi^0 = 1$ on a neighborhood of ${\rm supp}(\psi^0)$ and so that
$$
{\rm supp}(\tilde\psi^0)\subset\{|\xi|<2\}\cup\{|\xi_{2n-1}|< \frac{3}{4}| \xi'|\}.
$$ 
The operator $\Psi$ is defined by
$$
\widehat{\Psi u}(\xi) = \psi (\xi )\hat{u}(\xi).
$$
The operators $\Psi^+,\ \Psi^-,\ \Psi^0,\ \tilde\Psi^+,\ \tilde\Psi^-$, and $\tilde\Psi^0$ are defined
as above with  substitution of $\psi^+,\ \psi^-,\ \psi^0,\ \tilde\psi^+,\ \tilde\psi^-$, and $\tilde\psi^0$
for $\psi$, respectively.

The microlocal decomposition $\varphi=\varphi^++\varphi^-+\varphi^0$, alluded to in the
introduction, is now interpreted as follows:
$$
\varphi =\zeta\Psi^+\varphi +\zeta\Psi^-\varphi +\zeta\Psi^0\varphi,
$$
for all $\varphi\in C_0^\infty (U)$, where $\zeta\in C^\infty_0(U')$, $\bar{U}\subset U'$ and $\zeta=1$ on $U$. The
micro\-localization
$\varphi^+$ is also sometimes interpreted as
$\zeta\tilde\Psi^+\varphi$ and similarly with\break $\varphi^-$ and $\varphi^0$.

The following lemma is a consequence of the general G\aa rding inequality; here we give a proof which does not invoke the
general case.

\proclaim{Lemma} Let $(a_{IJ})$ be a matrix of $C^\infty$ functions on $U'$ which is non{\rm -}negative. Let 
$U\subset\bar U\subset U'$, $\zeta\in C_0^\infty(U')$ with $\zeta=1$
on $U$ and $\sigma\in C^\infty (U')$ with $\sigma\le 2$. Then{\rm ,} if
$U'$ is sufficiently small{\rm ,}  \pagebreak
\begin{eqnarray*}
\sum_{IJ}(a_{IJ}\sigma T\zeta\Psi^+\varphi_I,\sigma\zeta\Psi^+\varphi_J)
&\ge&-C\Bigl(\max_{{\rm supp}\, \zeta}|D(\sigma)|+1)\|\zeta\Psi^+\varphi\|^2\\
&&+\ C_\sigma\|\Psi^+\varphi\|_{-1}\|\zeta\Psi^+\varphi\|+\|\Psi^+\varphi\|_{-1}^2\Bigr),
\end{eqnarray*}
and
\begin{eqnarray*}
\sum_{IJ}(a_{IJ}\sigma T\zeta\Psi^-\varphi_I,\sigma\zeta\Psi^-\varphi_J)
&\ge&-C\Bigl(\max_{{\rm supp}\, \zeta}|D(\sigma)|+1)\|\zeta\Psi^-\varphi\|^2\\
&&+\ C_\sigma\|\Psi^-\varphi\|_{-1}\|\zeta\Psi^-\varphi\|+\|\Psi^-\varphi\|_{-1}^2\Bigr),
\end{eqnarray*}
where 
\begin{equation}
C_\sigma=\max_{{\rm supp}\, \zeta}|D^2(\sigma)|+\max_{{\rm supp}\, \zeta}|D(\sigma)|^2+1,
\end{equation}
for all $\varphi$ with ${\rm supp}(\varphi )\subset U$. Here $D$ denotes first
partial derivatives.
\endproclaim

{\it Proof}.  Let $\theta$ be a conical cutoff function with ${\rm supp}\ \theta\subset {\rm supp}\ \tilde\psi^+$ and $\theta=1$ on ${\rm supp}\
\psi^+$; denoting by $\Theta$ the corresponding pseudodifferential operator, we have
$$
\zeta\Psi^+\varphi=\zeta(\tilde\Psi^+)^2\Theta\Psi^+\varphi=(\tilde\Psi^+)^2\zeta\Psi^+\varphi+
[\zeta,(\tilde\Psi^+)^2]\Theta\Psi^+\varphi.
$$
Then, since the supports of the symbols of
$\Theta$ and of $[\zeta ,(\tilde{\Psi}^+)^2]$ are disjoint, the operator $ [\zeta ,(\tilde{\Psi}^+)^2]\Theta$
is of order $-\infty$ and we have  
\begin{eqnarray*}
\sum_{IJ}(a_{IJ}\sigma T\zeta\Psi^+\varphi_I,\sigma\zeta\Psi^+\varphi_J)&=&
\sum_{IJ}(a_{IJ}\sigma \tilde\zeta T(\tilde{\Psi}^+)^2\zeta\Psi^+\varphi_I,\sigma\zeta\Psi^+\varphi_J)\\
&&+\ O(\|\Psi^+\varphi\|^2_{-1}),
\end{eqnarray*}
where $\tilde\zeta\in C_0^\infty(U')$ with $\tilde\zeta=1$ on ${\rm supp}\ \zeta$. Let $R$ denote the pseudodifferential operator
of order
$1/2$ whose symbol is
$\xi_{2n-1}^{\frac{1}{2}}\psi^+(\xi)$;  then $T(\tilde{\Psi}^+)^2=R^*R$ and 
\begin{eqnarray*}
\sum_{IJ}(a_{IJ}\sigma T\zeta\Psi^+\varphi_I,\sigma\zeta\Psi^+\varphi_J)&=&
\sum_{IJ}\Bigl((a_{IJ}\sigma \tilde\zeta R\zeta\Psi^+\varphi_I,\sigma\tilde\zeta R\zeta\Psi^+\varphi_J)\\
&&+\ ([a_{IJ}\sigma^2\tilde\zeta,R^*]R\zeta\psi^+\varphi_I,\zeta\psi^+\varphi_J)\\
&&+\ O(\|\Psi^+\varphi\|^2_{-1})\Bigr).
\end{eqnarray*}  
The first term on the right is nonnegative and from the
pseudodifferential operator calculus we obtain
$$
\|[a_{IJ}\sigma^2\tilde\zeta,R^*]R\zeta\psi^+\varphi_I\|\le
C\Bigl((\max_{{\rm supp}\,\zeta}|D(\sigma)|+1)\|\zeta\Psi^+\varphi\|+C_\sigma\|\zeta\Psi^+\varphi\|_{-1})\Bigr).
$$
The first inequality in the lemma then follows and the second is proved analogously. 
\pagebreak

\section{Weighted microlocal estimates}
 
We begin by discussing the estimate for $\zeta\Psi^+\varphi$ for $(0,1)$-forms; here the
calculation is more transparent because we do not have to deal with complicated indices. We will
then show how the proof generalizes to the case of $(p,q)$-forms. 

\proclaim{Lemma} Let $M$ be a pseudoconvex {\rm CR} manifold of dimension $2n-1${\rm ,} 
let $x_o\in M$ and let $\lambda$ be a real nonnegative $C^\infty$ function with $s\in\Bbb   R^+$ such that $s\lambda\le
1$ and ${\rm Re}\Bigl(L_i\bar{L}_j(\lambda)\Bigr)$ is positive\/{\rm -}\/definite. Then there exists a neighborhood $U$ of $x_o$ and
$C>0$ such that
\begin{eqnarray}
\qquad s\|\zeta\Psi^+\varphi\|^2 
&\hskip-6pt \le\hskip-6pt &
C\Bigl( Q_b(\zeta\Psi^+\varphi,\zeta\Psi^+\varphi)+(s^2\max_{{\rm supp}(\zeta)}D(\lambda)^2+1)
\|\zeta\Psi^+\varphi\|^2\\  
&\hskip-6pt\hskip-6pt&+\ C(s,\lambda)\|\Psi^+\varphi\|_{-1}\|\zeta\Psi^+\varphi\|
+\|\Psi^+\varphi\|^2_{-1}\Bigr),\nonumber
\end{eqnarray}
where 
$$C(s,\lambda)=s\max_{{\rm supp}(\zeta)}|D^2(\lambda)|+s^2\max_{{\rm supp}(\zeta)}D(\lambda)^2+1,$$
$\zeta\in C^\infty_0(U)${\rm ,} $\varphi\in {\cal A}^{0,1}_b${\rm ,} ${\rm supp}(\varphi )\subset U'${\rm ,} and
$U'\supset\bar{U}$. 
\endproclaim 

 \vglue9pt
{\it Proof}.  If $\varphi\in C^\infty_0(U)$
then: 
\begin{eqnarray}
\|\sigma\bar{\partial}_b\varphi\|^2+\|\sigma\bar{\partial}_b^*\varphi\|^2&=&
\sum_{i<j}\|\sigma(\bar{L}_j\varphi_i -\bar{L}_i\varphi_j)\|^2
\\ && + \ \|\sigma\sum\bar{L}^*_i{\varphi}_i\|^2\nonumber\\
&&+\ O\lrp{\sum\|\sigma\bar{L}_i{\varphi}_j\|\|\sigma\varphi\|+\|\sigma\varphi\|^2}.\nonumber
\end{eqnarray}
We have $\bar{L}_i^*=-L_i+a_i$. Substituting this in the above we use the following integrations by
parts to ``convert'' the $L$ into $\bar{L}$:
$$
(\sigma^2L_iu,v)=-(\sigma^2u,\bar{L}_iv)+O(\|L_i(\sigma)u\|\|\sigma v\|+|\sigma u\|\|\sigma v\|)
$$
and
\begin{eqnarray*}
(\sigma^2\bar L^*_iu,\bar L^*_jv)&=&\bigg(\sigma^2L_iu,L_jv)+O(\|\sigma u\|\|\sigma v\|\\
&  &+\ \sum(\|\sigma\bar L_ku\|\|\sigma v\|+\|\sigma u\|\|\sigma\bar L_kv\|\bigg)\\
&=&(\sigma\bar{L}_ju,\sigma\bar{L}_iv)+(\sigma c_{ij}Tu,\sigma v)+
2(L_i\bar{L}_j(\sigma)u,\sigma v)\\&&+\ O\Bigl(\|\sigma u\|\|\sigma v\|+\|\sigma\bar L_ju\|\|L_i(\sigma) v\| \\
&&+\ \|L_i(\sigma)u\|\|\sigma\bar L_jv\|+\sum(\|\sigma\bar L_ku\|\|\sigma v\|+\|\sigma u\|\|\sigma\bar L_kv\|)\Bigr).
\end{eqnarray*}
\noindent
Let $\sigma_s=e^{s\lambda}$; then, substituting $\sigma_s$ for $\sigma$, we have \pagebreak

\begin{eqnarray*}  \noalign{\vskip-16pt}
\|\sigma_s\bar{\partial}_b^*\varphi\|^2&=&\|\sum_i\sigma_s\bar L_i^*\varphi_i\|^2\\&=&
\sum_i\|\sigma_s\bar L_i^*\varphi_i\|^2+2\sum_{i<j}{\rm Re}(\sigma_s^2\bar L_i^*\varphi_i,\bar L_j^*\varphi_j)\\
&=&\sum_i\|\sigma_s\bar L_i\varphi_i\|^2+ 2\sum_{i<j}{\rm Re}(\sigma_s^2\bar L_j\varphi_i,\bar L_i\varphi_j)\\
&&+\ (\sigma_s c_{ij}T\varphi_i,\sigma_s\varphi_j)+(L_i\bar{L}_j(\sigma_s)\varphi_i,\sigma_s\varphi_j)\\
&&+\ O\Bigl(\|\sigma_s\varphi\|^2+\sum_{i,j}\|\sigma_s\bar L_j\varphi_i\|\|L_i(\sigma_s)\varphi_j\|+
\sum_k(\|\sigma_s\bar L_k\varphi\|\|\sigma_s\varphi\|)\Bigr).
\end{eqnarray*}
Combining this with (3.2) we obtain
\begin{eqnarray*}
\|\sigma_s\bar{\partial}_b\varphi\|^2+\|\sigma_s\bar{\partial}_b^*\varphi\|^2&=&
\sum_{i,j}
\Bigl(\|\sigma_s\bar L_i\varphi_j\|^2+(\sigma_s
c_{ij}T\varphi_i,\sigma_s\varphi_j)\\
&&+\ {\rm Re}(L_i\bar{L}_j(\sigma_s)\varphi_i,\sigma_s\varphi_j)\Bigr)\nonumber\\
&&+\ O\Bigl(\|\sigma_s\varphi\|^2+\sum_{i,j}\|\sigma_s\bar L_j\varphi_i\|\|L_i(\sigma_s)\varphi_j\|\nonumber\\
&&+\ \sum_k(\|\sigma_s\bar L_k\varphi\|\|\sigma_s\varphi\|)\Bigr).
\end{eqnarray*}
Since $\sigma_s$ is bounded independently of $s$,
\begin{eqnarray}
&&\hskip-.5in
\sum_{i,j}\Bigl((\sigma_sc_{ij}T\varphi_i,\sigma_s\varphi_j)+{\rm Re}(L_i\bar{L}_j(\sigma_s)\varphi_i,\sigma_s\varphi_j)\Bigr)
\\
&&\qquad \le C\Bigl(Q_b(\varphi,\varphi) 
+\|\varphi\|^2+\sum_{i,j}\|L_i (\sigma_s)\varphi_j\|^2\Bigr).\nonumber \end{eqnarray}
Note that $D(\sigma_s)=sD(\lambda)\sigma_s$ and $D^2(\sigma_s)=sD^2(\lambda)\sigma_s+s^2D(\lambda)^2\sigma_s$. Then we have
\begin{eqnarray*}
\max_{{\rm supp}(\zeta)}|D(\sigma_s)|&\sim& s\max_{{\rm supp}(\zeta)}|D(\lambda)|,\\
 \max_{{\rm supp}(\zeta)}|D^2(\sigma_s)|&\sim& s\max_{{\rm supp}(\zeta)}|D^2(\lambda)|+s^2\max_{{\rm
supp}(\zeta)}D(\lambda)^2,\\
 {\rm Re}  L_i\bar L_j(\sigma_s)&\sim& s{\rm Re}   L_i\bar L_j(\lambda)+s^2D(\lambda)^2,\end{eqnarray*}
and
$$C_{\sigma_s}\sim C(s,\lambda)=s\max_{{\rm supp}(\zeta)}|D^2(\lambda)|+s^2\max_{{\rm supp}(\zeta)}D(\lambda)^2+1,$$
where $C_{\sigma_s}$ is defined by (2.1).

When we substitute $\zeta\Psi^+\varphi$ for $\varphi$ in (3.3) the first term is estimated from below by\pagebreak
\begin{eqnarray*}
&&\hskip-.25in \sum_{i,j}(\sigma_sc_{ij}T\zeta\Psi^+\varphi_i,\sigma_s\zeta\Psi^+\varphi_j) \\
&&\ge  
-C\Bigl( (s\max_{{\rm supp}\zeta}|D(\lambda)|\|\zeta\Psi^+\varphi\|^2  
 +\ C(s,\lambda)\|\Psi^+\varphi\|_{-1}\|\zeta\Psi^+\varphi\|
+\|\Psi^+\varphi\|^2_{-1}\Bigr).
\end{eqnarray*}
The second term is estimated from below by
$$(L_i\bar{L}_j(\sigma_s)\zeta\Psi^+\varphi_i,\sigma_s\zeta\Psi^+\varphi_j)\ge sC'\|\sigma_s\zeta\Psi^+\varphi\|^2
-s^2\max_{{\rm supp}\zeta}D(\lambda)^2\|\zeta\Psi^+\varphi\|^2.$$
Since $\sigma_s$ is bounded away from zero independently of $s$ the estimate (3.1) follows from the above and (3.3) thus
proving the lemma.\hfill\qed
\vglue9pt

To generalize the above to $(0,q)$-forms $\varphi=\sum\varphi_I\bar\omega_I$ we define
$$
A_{IJ}(\lambda)={\rm Re}\sum_{i,j,K}\epsilon^{iK}_I\epsilon^{jK}_JL_i\bar{L}_j(\lambda)
$$
and
$$
c_{IJ}=\sum_{i,j,K}\epsilon^{iK}_I\epsilon^{jK}_Jc_{ij},
$$
where $K$ runs over all ordered $(q-1)$-tuples. Each of the coefficients $\epsilon^{iK}_I$ is either
$0, 1,$ or $-1$ defined as follows. First, if $i\notin K$ we denote by $<iK>$ the ordered
$q$-tuple containing $i$ and the elements of $K$. Then we define
$$
\epsilon^{iK}_I=\left\{ \begin{array}{ll}  0 &\qquad \hbox{if $i\in K$} \\ 0 &\qquad \hbox{if $\langle iK\rangle \neq I$}\\
\mathop{{\rm sgn}}\nolimits\left({iK\atop I}\right) &\qquad \hbox{if $\langle iK\rangle =I$}, \end{array}\right.
$$
where $\mathop{{\rm sgn}}\nolimits\left({iK\atop I}\right)$ denotes the sign of the permutation
$\{i,K\}\to I$.

\proclaim{Lemma} \hskip-6pt Let $M$ be a pseudoconvex {\rm CR} manifold of dimension $2n-1${\rm ,} let $x_o\in M$ and
let
$\lambda$ be a real nonnegative $C^\infty$ function{\rm ,} $s\in\Bbb   R^+$ with $s\lambda\le 1${\rm ,} and suppose that
${\rm Re} A_{IJ}(\lambda)$ is positive\/{\rm -}\/definite. Then there exists a neighborhood $U$ of $x_o$ and
$C>0$ such that
\begin{eqnarray}
\qquad s\|\zeta\Psi^+\varphi\|^2 
& \le &
C\Bigl(Q_b(\zeta\Psi^+\varphi,\zeta\Psi^+\varphi)+(s^2\max_{{\rm supp}(\zeta)}
D(\lambda)^2+1)\|\zeta\Psi^+\varphi\|^2\\  
&&+\ C(s,\lambda)\|\Psi^+\varphi\|_{-1}\|\zeta\Psi^+\varphi\|
+\|\Psi^+\varphi\|^2_{-1}\Bigr),\nonumber
\end{eqnarray}
where 
$$C(s,\lambda)=s\max_{{\rm supp}(\zeta)}|D^2(\lambda)|+s^2\max_{{\rm supp}(\zeta)}D(\lambda)^2+1,$$
$\zeta\in C^\infty_0(U)$, $\varphi\in {\cal A}^{0,q}_b$, ${\rm supp}(\varphi )\subset U'$, and $U'\supset\bar{U}$.  
\endproclaim

The proof of this lemma is entirely analogous to the case of $(0,1)$-forms. Integrating by parts as
above we obtain the following in place of (3.2) 
\begin{eqnarray}
&&{\rm Re}
\sum_{I,J}\Bigl((\sigma_sc_{IJ}T\varphi_I,\sigma_s\varphi_J)+(A_{IJ}(\sigma_s)\varphi_I,\sigma_s\varphi_J)\Bigr)
\\
&&\qquad \le C\Bigl(Q_b(\varphi,\varphi)
 +\ \|\varphi\|^2+\|D(\sigma_s)\varphi\|^2\Bigr).\nonumber
\end{eqnarray}
Observe that $c_{ij}\ge 0$ implies that $c_{IJ}\ge 0$ since at a point $x_o$ we can choose the
$\{L_i\}$ so that $c_{ij}(x_o)=\delta_{ij}c_{ii}(x_0)$ and then $c_{IJ}(x_o)=\delta_{IJ}\sum_{i\in
I}c_{ii}(x_o)$. Hence, again substituting $\zeta\Psi^+\varphi$ for $\varphi$, we can apply Lemma 2.5
and conclude the proof of Lemma 3.2. \hfill\qed\vglue12pt

To microlocalize with $\Psi^-$ we define the local conjugate-linear duality map
$F^q:{\cal A}^{0,q}_b\to{\cal A}^{0,n-q-1}_b$ as follows. If
$\varphi=\sum\varphi_I\bar{\omega}_I$ then
\begin{equation}
F^q(\varphi)=\sum\epsilon^{II'}\bar{\varphi}_{I}\bar{\omega}_{I'},
\end{equation}
where $I'$ denotes the strictly increasing $(n-q-1)$-tuple consisting of all integers in $[1,n-1]$
which do not belong to $I$ and $\epsilon^{II'}$ is the sign of the permutation $\{I,I'\}\to \langle \{1,\dots, n-1\}\rangle $.
Then
$F^{n-q-1}F^q\varphi=\varphi$ and we have
\begin{eqnarray*}
\bar{\partial}_bF^q\varphi &=&F^{q-1}\bar{\partial}^\ast_b\varphi +\sum a_{IJ}\bar\varphi_I\bar{\omega}_J
\\
\noalign{\noindent 
and}
\bar{\partial}^\ast_b F^q\varphi &=&F^{q+1}\bar{\partial}_b\varphi +\sum b_{HK}\bar\varphi_H\bar{\omega}_K.
\end{eqnarray*}
Hence 
$$
Q_b(F^q\varphi,F^q\varphi)=O(Q_b(\varphi ,\varphi)+\|\varphi\|^2).
$$
Thus replacing $\varphi$ by $F^q\varphi$ in (3.5), we obtain
\begin{eqnarray}
&&\hskip-.25in {\rm Re}\sum_{I,J}\Bigl((\sigma_sc_{I'J'}T\bar\varphi_I,\sigma_s\bar\varphi_J)
+(A_{I'J'}(\sigma_s)\bar\varphi_I,\sigma_s\bar\varphi_J)\Bigr)
\\ &&\qquad \le C\Bigl(Q_b(\varphi,\varphi)+\|\varphi\|^2+\|D(\sigma_s)\varphi\|^2\Bigr).\nonumber
\end{eqnarray} 
Now note that, since $\bar{T}=-T$, 
$$ 
{\rm Re} \sum_{I,J}(\sigma_sc_{I'J'}T\bar{\varphi}_I,\sigma_s\bar{\varphi}_J)=
-{\rm Re} \sum_{I,J}(\sigma_sc_{J'I'}T\varphi_I,\sigma_s\varphi_J).
$$
Hence, substituting $\zeta\Psi^-\varphi$ for $\varphi$, we proceed as above and obtain the
following result.\pagebreak

\proclaim{Lemma} Let $M$ be a pseudoconvex {\rm CR} manifold of dimension\break $2n-1${\rm ,} let $x_o\in M$ and let
$\lambda$ be a real non-negative $C^\infty$ function, $s\in\Bbb   R^+$ with $s\lambda\le 1$, and suppose that
${\rm Re} A_{I'J'}(\lambda)$ is positive\/{\rm -}\/definite. Then there exists a neighborhood $U$ of $x_o$ and
$C>0$ such that
\begin{eqnarray}
\qquad s\|\zeta\Psi^-\varphi\|^2 
& \le &
C\Bigl(Q_b(\zeta\Psi^-\varphi,\zeta\Psi^-\varphi)+(s^2\max_{{\rm supp}(\zeta)}D(\lambda)^2+1)\|\zeta\Psi^-\varphi\|^2
\\  
&&+\ C(s,\lambda)\|\Psi^-\varphi\|_{-1}\|\zeta\Psi^-\varphi\|
+\|\Psi^-\varphi\|^2_{-1}\Bigr),\nonumber
\end{eqnarray}
where 
$$C(s,\lambda)=s\max_{{\rm supp}(\zeta)}|D^2(\lambda)|+s^2\max_{{\rm supp}(\zeta)}D(\lambda)^2+1,$$
$\zeta\in C^\infty_0(U)${\rm ,} $\varphi\in {\cal A}^{0,q}_b${\rm ,} ${\rm supp}(\varphi )\subset U'${\rm ,} and $U'\supset\bar{U}$.  
\endproclaim

\section{The localization lemma}
 
{\it Definition} 4.1.
Let $S\subset M$ be a submanifold of $M$; the {\it holomorphic dimension} of $S$ at
$x_o\in S$ is the dimension of $T^{1,0}_{x_o}(M)\cap{\Bbb C}{\rm T}_{x_o}(S)$ and  $S$ is called
{\it totally real} if the holomorphic dimension of $S$ at $x$ is zero for all $x\in S$. 
\vglue4pt\advance\theoremcount by 1

We are now in a position to formulate Lemma 1.4 precisely.

\proclaimtitle{Localization Lemma}
\proclaim{Lemma}
Let $M$ be a pseudoconvex {\rm CR} manifold and let $S\subset M$ be a submanifold of holomorphic dimension less
than or equal to $m$. Let $x_0\in S$. For small positive $a$ let $S_a$ denote the set of $x\in M$
such that $d(x,S)<a${\rm ,} where $d(x,S)$ denotes the distance from $x$ to $S$. Then there exist positive
constants $a_o,\  C$ and neighborhoods $U${\rm ,} $U'$ of $x_0${\rm ,} with $\bar{U}\subset U'$ and a constant $C$ such that 
\begin{equation}
\|\zeta\Psi^+ \varphi \|^2_{S_a} \le C\Bigl( a^2Q_b(\zeta\Psi^+\varphi, \zeta\Psi^+\varphi) + 
\| \zeta\Psi^+\varphi \|^2_{M -S_a}+\|\Psi^+\varphi\|^2_{-1}\Bigr),\quad 
\end{equation}
for $0<a<a_o${\rm ,} $\zeta\in C^\infty_0(U)$ and all $(p,q)$\/{\rm -}\/forms $\varphi$ with $n-1\ge q\ge m+1$
supported in
$U'$. Furthermore
\begin{equation}
\| \zeta\Psi^-\varphi \|^2_{S_a} \le C\Bigl( a^2Q_b( \zeta\Psi^-\varphi,  \zeta\Psi^-\varphi) + 
\|  \zeta\Psi^-\varphi \|^2_{M -S_a}+\|\Psi^-\varphi\|^2_{-1}\Bigr),\quad 
\end{equation}
for all $(p,q)$\/{\rm -}\/forms $\varphi$ with $0\le q\le n-m-2$ supported in $U'$. Here $\|\ \|_{S_a}$\ and
$\|\ \|_{M-S_a}$ denote $L_2$\/{\rm -}\/norms over $S_a$ and $M-S_a$\/{\rm ,} respectively.    
\endproclaim

\demo{Proof}    Let $f_1,\ldots ,f_k$ be real-valued functions in a neighborhood of $x_o$ such that $x\in S$ if and only
if $f_i(x)=0$ for $i=1,\ldots ,k$. Suppose further that the gradients of the $f_i$ are linearly
independent. Since the rank of $(L_if_j|_{x_o})$ is greater than or equal to $n-1-m$, we can choose 
( after renumbering the $L_i$) functions $\{g_1,\ldots g_{n-1-m}\}$ which are real linear combinations of
the
$f_i$, such that  $L_i(g_j)|_{x_o}=\delta_{ij}$ when $i =1,\ldots n-q-m$. Then, setting
$\lambda =\sum g_k^2$ we have $A_{IJ}(\lambda )|_{x_o}=2\delta_{IJ}p_I$, where $I$ and $J$ are
ordered $q$-tuples so that $p_I$ is greater than or equal to the number of elements in
$I\cap\{1,\ldots ,n-1-m\}$. Since $m\le q-1$ we have $p_I>0$ and hence  $(A_{IJ}(\lambda )|_{x_o})$ is
positive definite. Thus there exists a neighborhood $U'$ of $x_o$ on which ${\rm Re} (A_{IJ}(\lambda))$
is positive definite.

Let $\zeta_a\in C_0^\infty(S_{2a})$ with $\zeta_a=1$ on $S_a$ such that
$|D\zeta_a|\le\frac{C}{a}$. Then  $\lambda\le Ca$ and $|L_i(\lambda)|\le Ca$ in $S_{2a}$. Setting
$s=\frac{\delta}{a^2}$ with $\delta$ sufficiently small so that $s\lambda\le 1$ in $S_{2a}$ we can
apply Lemma 3.2 and obtain
\begin{eqnarray}
s\|\zeta_a\zeta\Psi^+\varphi\|^2 
& \le &
C\Bigl(Q_b(\zeta_a\zeta\Psi^+\varphi,\zeta_a\zeta\Psi^+\varphi)\\
&&+
\ \Bigl(s^2\max_{{\rm supp}(\zeta_a\zeta)}D(\lambda)^2+1)\|\zeta_a\zeta\Psi^+\varphi\|^2\nonumber\\  
&&+\ C(s,\lambda)\|\Psi^+\varphi\|_{-1}\|\zeta_a\zeta\Psi^+\varphi\|
+\|\Psi^+\varphi\|^2_{-1}\Bigr),\nonumber
\end{eqnarray}
where 
$$C(s,\lambda)=s\max_{{\rm supp}(\zeta_a\zeta)}|D^2(\lambda)|+s^2\max_{{\rm supp}(\zeta_a\zeta)}D(\lambda)^2+1,$$
$\zeta\in C^\infty_0(U)$, $\varphi\in {\cal A}^{0,q}_b$, ${\rm supp}(\varphi )\subset U'$, and $U'\supset\bar{U}$.
Dividing both sides of (4.3) by $s$ and  with $s=\frac{\delta}{a^2}$ and $a_o\le\delta$ we have
$$s\max_{{\rm supp}(\zeta_a\zeta)}D(\lambda)^2+\frac{1}{s}\le C\delta$$
and 
$$\frac{1}{s}C(s,\lambda)\le C.$$
Hence
\begin{eqnarray*}
\|\zeta_a\zeta\Psi^+\varphi\|^2 
& \le &
C\Bigl(\frac{a^2}{\delta}Q_b(\zeta_a\zeta\Psi^+\varphi,\zeta_a\zeta\Psi^+\varphi)+
\delta\|\zeta_a\zeta\Psi^+\varphi\|^2\\  
&&+\ \|\Psi^+\varphi\|_{-1}\|\zeta_a\zeta\Psi^+\varphi\|
+\|\Psi^+\varphi\|^2_{-1}\Bigr). 
\end{eqnarray*}
Thus, since
$$\|\zeta\Psi^+\varphi\|_{S_a}\le\|\zeta_a\zeta\Psi^+\varphi\|$$
and
$$a^2Q_b(\zeta_a\zeta\Psi^+\varphi,\zeta_a\zeta\Psi^+\varphi)\le a^2Q_b(\zeta\Psi^+\varphi,\zeta\Psi^+\varphi)+
C\|\zeta\Psi^+\varphi\|_{M-S_a}$$
then (4.1) follows by choosing $\delta$ sufficiently small.

To prove (4.2) note that $A_{I'J'}(\lambda )|_{x_o}=2\delta_{I'J'}p_{I'}$, where $I'$ and $J'$ are
ordered $(n-1-q)$-tuples so that $p_{I'}$ is greater than or equal to the number of elements in
$I'\cap\{1,\ldots ,n-1-m\}$. Since $m\le n-q-2$ we have $p_{I'}
>0$ and hence  $(A_{I'J'}(\lambda
)|_{x_o})$ is positive definite. Then, proceeding as above, using Lemma 3.4 instead of Lemma 3.3 we obtain
the desired result. \enddemo

Note that for $\zeta\Psi^0\varphi$ the above estimates are a consequence of ellipticity and thus hold for all $q$. Then, when
the above estimates hold for $\zeta\Psi^+\varphi$ and $\zeta\Psi^-\varphi$, we obtain the following:
$$
\| \varphi \|^2_{S_a} \le C\Bigl( a^2Q_b(\varphi,\varphi) + 
\| \varphi \|^2_{M -S_a}+\|\varphi\|^2_{-1}\Bigr).
$$
Substituting $2a$ for $a$ and $\zeta_a\varphi$ for $\varphi$ we get
$$
\|\zeta_a \varphi \|^2 \le C\Bigl( a^2Q_b(\zeta_a\varphi,\zeta_a\varphi) +\|\zeta_a\varphi\|^2_{-1}\Bigr).
$$
Choosing $U'$ with sufficiently small diameter we get 
$$
C\|\zeta_a\varphi\|^2_{-1}\le\frac{1}{2}\|\zeta_a\varphi\|^2.
$$
Thus we obtain:

\proclaim{{C}orollary}
Let $M$ be pseudoconvex {\rm CR} manifold and let $S\subset M$ be a submanifold of holomorphic dimension less
than or equal to $m$ with $m\le \frac{1}{2}(n-3)$. Let $x_0\in S$. For small positive $a$ let $S_a$ denote the set of $x\in M$
such that $d(x,S)<a${\rm ,}
 where $d(x,S)$ denotes the distance from $x$ to $S$. Then there exist positive
constants $a_o,\  C,$ and a neighborhood $U$ such that 
\begin{equation}
\| \varphi \|^2_{S_a} \le C\Bigl( a^2Q_b(\varphi,\varphi) + 
\| \varphi \|^2_{M -S_a}\Bigr) 
\end{equation}
for $0<a<a_o$ and all $(p,q)$\/{\rm -}\/forms $\varphi$ with $m+1\le q\le n-m-2$
supported in~$U$.
\endproclaim

M. Christ has pointed out that pseudoconvexity is a crucial assumption for the estimate (4.4). He has supplied the
following example. Let $M$ be a five-dimensional CR manifold with coordinates $(z_1,z_2,t)$, whose CR\  structure is
defined by the vector fields
$$
 L_1 = \partial_{ z_1} -i\bar z_1 \partial_t, 
\qquad
 L_2 = \partial_{z_2} + i\bar z_2 \partial_t. 
$$
Then we have $\bar L_j^*=-L_j$; setting $T=i\partial_t$ we have $[L_1,\bar L_1]=2T$ and $[L_2,\bar L_2]=-2T$. Let $S$
be defined by $z_1=z_2=0$; then $m$, the holomorphic dimension of $S$, equals zero. Let $S_a=\{(z_1,z_2,t)\Big |\  |z|<
a\}$, with $|z|^2=|z_1|^2+|z_2|^2$. Let $g=|z|^2-it$; then $L_1g=\bar L_2g=0$. For each $\tau>0$ we define the
$(0,1)$-form $\varphi^\tau$ by $\varphi^\tau=f^\tau\zeta\bar\omega _1$, where $f^\tau (z_1,z_2,t) = \exp\tau(-g+g^2)$ and
$\zeta\in C_0^\infty(|z|^2+t^2<2r^2)$ with $\zeta=1$ on $|z|^2+t^2\le r^2$. Note that $L_1f^\tau =\bar L_2f^\tau=0$. We have
$$
Q_b(\varphi^\tau,\varphi^\tau)=\|f^\tau L_1\zeta\|^2+\|f^\tau\bar L_2\zeta\|^2.
$$
If $a<r\le\frac{1}{2}$ we have 

$$
{\rm Re}(-g+g^2)=-|z|^2-t^2+|z|^4\ge-\frac{a^2}{2}
$$
on the region $|z|^2+t^2\le\frac{a^2}{2}$ which is contained in $S_a$ so that
$$
ca^5\exp\lrp{-\frac{1}{2}\tau a^2}\le\|\varphi^\tau\|^2_{S_a}.
$$
Furthermore
$$
{\rm Re}(-g+g^2)=-|z|^2-t^2+|z|^4\le-\frac{3}{4}r^2\le-\frac{3}{4}a^2
$$
on the support of the derivatives of $\zeta$; thus
$$
Q_b(\varphi^\tau,\varphi^\tau)+\|\varphi^\tau\|^2_{M-{S_a}}\le C\exp\lrp{-\frac{3}{4}\tau a^2}.
$$
This contradicts (4.3) for large $\tau$.

\section{Subelliptic multipliers}

To study the operators $\Box_b$, $\bar{\partial}_b$, and the $\bar{\partial}$-Neumann  problem 
on $(p,q)$-forms we define microlocal subelliptic multipliers below. In this section we develop 
those properties of the multipliers which are needed in the proofs of superlogarithmic
estimates and the hypoellipticity results. Explicit formulas for subelliptic multipliers are given in the appendix.  
\numbereddemo{Definition}
For $x_o\in M$ let ${\cal I}^+_q(x_o)$ be the subset of germs of $C^\infty$ functions at $x_o$
defined as follows. The germ $\rho$ is in ${\cal I}^+_q(x_o)$ if and only if there exist $U, U',
\varepsilon, C$, and $\zeta$ such that
\begin{eqnarray}
\| \rho \zeta \Psi^+ \varphi\|^2_\varepsilon \le C\bigl ( Q_b(\zeta \Psi^+ \varphi, \zeta \Psi^+
\varphi) + \|\Psi^+\varphi \|^2\bigr )
\end{eqnarray}
for all $(0,q)$-forms $\varphi$ with ${\rm supp}(\varphi)\subset U'$. Here, $U$ is a neighborhood of $x_o$
such that $\bar{U}\subset U'$, $\varepsilon$ and $C$ are positive constants and $\zeta \in
C_0^\infty(U)$ such that $\zeta =1$ on a neighborhood of $x_o$. We define ${\cal I}^0_q(x_o)$ and
${\cal I}^-_q(x_o)$ by replacing
$\Psi^+$ in the above by $\Psi^0$ and $\Psi^-$, respectively. 
\enddemo
Note that ${\cal I}^0_q(x_o)$ is the set of  all germs of $C^\infty $ functions at $x_o$; in fact
we have
$$
\| \rho \zeta \Psi^0 \varphi\|^2_1 \le C\bigl ( Q_b(\zeta \Psi^0 \varphi, \zeta \Psi^0
\varphi) + \|\Psi^0\varphi \|^2\bigr ).
$$
Thus, if we denote the set of subelliptic multipliers satisfying (5.1) by ${\cal I}_q(x_o)$, we
have ${\cal I}_q(x_o) = {\cal I}^+_q(x_o) \cap {\cal I}^-_q(x_o)$.

 \proclaim{Proposition}
  If $1\le q' \le q \le n-1$ then ${\cal I}^+_q(x_o) \supset
{\cal I}^+_{q'}(x_o)$.
\endproclaim

{\it Proof}. Now, $\rho \in {\cal I}^+_{q-1}(x_o)$ and we want to show that $\rho \in
{\cal I}^+_q(x_o)$. Let $\varphi$ be a $(0,q)$-form supported in $U$,
$$
\varphi = \sum {\varphi}_I\bar{\omega}_I,
$$
where the $I$ are ordered $q$-tuples of integers between $1$ and $n-1$. We define the
$(0,q-1)$-forms $\beta^j$ by
$$
\beta^j =\sum_K \varphi_{<Kj>}\bar{\omega}_K,
$$
where $I=\langle Kj\rangle $ and $q\le j \le n-1$. Thus we have

$$
\| \rho \zeta \Psi^+ \beta^j\|^2_\varepsilon \le C\bigl ( Q_b(\zeta \Psi^+ \beta^j, \zeta \Psi^+
\beta^j) + \| \varphi \|^2\bigr ).
$$
\noindent
Now

$$  Q_b(\zeta \Psi^+ \beta^j, \zeta \Psi^+\beta^j) =
\sum \|\bar{L}_k\Psi^+\beta^j_K\|^2 + \sum(b_{HK}T\zeta\Psi^+\beta^j_H,\zeta\Psi^+\beta^j_K) +
O(\|\varphi\|^2),
$$
\noindent
where $b_{HK}=\sum\epsilon^{kH}_{hK}c_{hk}$ and
$$  Q_b(\zeta \Psi^+ \varphi, \zeta \Psi^+\varphi) =
\sum \|\bar{L}_k\Psi^+\varphi_I\|^2 + \sum(c_{IJ}T\zeta\Psi^+\varphi_I,\zeta\Psi^+\varphi_J) +
O(\|\varphi\|^2),
$$
where $c_{IJ}=\sum\epsilon^{jI}_{iJ}c_{ij}$. To compare $a_{IJ}$ and $b_{HK}$ at a point $x$ we
choose the $L_1,\ldots ,L_{n-1}$ so that $c_{ij}(x)=\delta_{ij}c_{ii}(x)$.  Then we have
$$
c_{IJ}(x)=\delta_{IJ}\sum_{i\in I} c_{ii}(x) 
$$
and 
$$
b_{HK}(x)=\delta_{HK}\sum_{h\in H}c_{hh}(x).
$$
If $I=\langle Hk\rangle $,
$$
c_{II}(x_o)=b_{HH}(x_o)+c_{kk}(x_o)\ge b_{HH}(x_o);
$$
hence
$$
\sum_{H, K}b_{HK}(x)\beta^j_H(x)\bar{\beta}^j_K(x)\le
\sum_{IJ}c_{IJ}(x)\varphi_I(x)\bar{\varphi}_J(x).
$$
\noindent
Since
$$
\|\rho\zeta\Psi^+\varphi\|_\varepsilon \le \sum_j\| \rho\zeta\Psi^+\beta^j\|_\varepsilon
$$
we conclude (using G\aa rding's inequality) that
$$
\| \rho \zeta \Psi^+ \varphi\|^2_\varepsilon \le C\bigl ( Q_b(\zeta \Psi^+ \varphi, \zeta \Psi^+
\varphi) + \| \varphi \|^2\bigr )
$$
so that $\rho \in {\cal I}^+_q(x_o)$.  
\hfill\qed

\proclaim{Proposition}
 ${\cal I}^-_q(x_o)={\cal I}^+_{n-q-1}(x_o)$ for $0\le q\le n-2$.
\endproclaim
 
{\it Proof}. If $\varphi$ is a $(0,q)$-form with support in $U$ given by $\varphi =
\sum\epsilon^{II'}\varphi_I\bar{\omega}_I$, as in Section 3, we set $F_q\varphi =\sum
\bar{\varphi}_{I'}\bar{\omega}_{I'}$, where $I'$ denotes the ordered $(n-q-1)$-tuple consisting of
the complement of $I$. Then noting that $\| \Psi^+u\|=\| \Psi^-\bar{u}\|$, we have
$$
\| \rho\zeta\Psi^+\varphi\|_\varepsilon =\| \rho\zeta\Psi^-F_q\varphi\|_\varepsilon +
O(\|\varphi\|)
$$
and
$$Q_b(\zeta\Psi^+\varphi ,\zeta\Psi^+\varphi) = Q_b(\zeta\Psi^-F_q\varphi ,\zeta\Psi^-F_q\varphi)+
O(\| \varphi\|^2),
$$
which concludes the proof. \hfill\qed\vglue12pt

Combining these propositions we see that ${\cal I}^-_q(x_o) \subset {\cal I}^+_q(x_o)$, if 
$n-q-1\ge q$, that is if $q\le \frac{1}{2}(n-1)$. Hence ${\cal I}_q(x_o)={\cal I}^+_q(x_o)$
whenever $q\le \frac{1}{2}(n-1)$. In particular ${\cal I}_1(x_o)={\cal I}^+_1(x_o)$ when
$n\ge 2$.

\section{Superlogarithmic estimates}

In this section we will prove superlogarithmic estimates under the assumption that there exists a
subelliptic multiplier $\rho$ satisfying the condition 
$$
(*)\quad \quad \lim_{x\to S}d(x,S)\log\rho (x) = 0.
$$
We define the operator $\log\Lambda$ by
$$\widehat{\log\Lambda u}(\xi)=\frac{1}{2}(\log(1+|\xi|^2)\hat u(\xi)$$ 
and we have the following result.

\proclaim{Theorem} Let $M$ be a pseudoconvex {\rm CR} manifold. Let $S\subset M$ be a 
manifold such that the holomorphic dimension of $S$ at each point is less than or equal to $q-1$.
Let $x_o\in S${\rm ,} let $U'$ be a neighborhood of $x_o$ and suppose that $\rho \in C^\infty(U')$
with $\rho \in {\cal I}^+_q(x_o)$ satisfying $(*)$. Let $U$ be a neighborhood of $x_o$ with 
$\bar{U}\subset U'$ and $\zeta\in C^\infty_0(U)$. Then the  following superlogarithmic
estimate 
$({\rm SL}^+)_q$ holds. For each $\delta > 0$ there exist  $C_\delta${\rm ,}
such that
$$
({\rm SL}^+)_q\qquad\|(\log\Lambda)\zeta\Psi^+\varphi\|^2\le
\delta^2Q_b(\zeta\Psi^+\varphi,\zeta\Psi^+\varphi) +  C_\delta\|\Psi^+\varphi\|^2,
$$
for all  $\varphi\in {\cal A}^{0,q}_b$ with support in $U'$. Furthermore{\rm ,} if the
holomorphic dimension of $S$ at each point is less than or equal to
$n-q-2${\rm ,} and if  $\rho \in {\cal I}^-_q(x_o)$ satisfying  $(*)$  then the estimate 
$({\rm SL}^-)_q$ given by
$$
({\rm SL}^-)_q\qquad\|(\log\Lambda)\zeta\Psi^-\varphi\|^2\le
\delta^2Q_b(\zeta\Psi^-\varphi,\zeta\Psi^-\varphi) +  C_\delta\|\Psi^-\varphi\|^2,
$$
for all $\varphi\in {\cal A}^{0,q}$ with support in $U'$. 
\endproclaim

{\it Proof}. Let $\gamma_0, \tilde\gamma_0, \gamma, \tilde\gamma$ be nonnegative functions on
$C^\infty([0,\infty))$ such that ${\rm supp}(\gamma_0)\subset {\rm supp}(\tilde\gamma_0)\subset [0,2),\ 
{\rm supp}(\gamma)\subset {\rm supp}(\tilde\gamma)\subset [1,3),\  \tilde\gamma_0=1$ on ${\rm
supp}(\gamma_0),\break 
\tilde\gamma=1$ on ${\rm supp}(\gamma)$, and    
$$
\gamma_o(x)^2 + \sum_k \gamma (2^{-k}x)^2 = 1,
$$
when $x \ge 0$. When $k\ge 1$ we set $\gamma_k(x)=\gamma(2^{-k}x)$, 
and define the operators $\Gamma_k$ and $\widetilde{\Gamma}_k$ by
$$
\widehat{\Gamma_ku}(\xi) = \gamma_k(|\xi|)\hat{u}(\xi)
$$
and
$$
\widehat{\tilde{\Gamma}_ku}(\xi) = \widetilde{\gamma}_k(|\xi|)\hat{u}(\xi).
$$
Let $\rho$ be a subelliptic multiplier, satisfying $(*)$. Then for each $\delta
>0$ there exists
$A_\delta >0$ such that
$$
\frac{1}{\rho (x)} \le A_\delta \exp(\frac{\delta}{d(x,S)}).
$$
Then, using the fact that $|\xi| \sim 2^k$ when $\xi \in {\rm supp}(\gamma_k)$ and applying (5.1), we have
\begin{eqnarray*}
\| \Gamma_k\zeta\Psi^+\varphi\|_{M-S_a}& \le &
{\rm max}_{M-S_a}{\frac{1}{\rho}}\|\rho \Gamma_k\zeta(\Psi^+\varphi)\|
\\& \le &A_\delta \exp(\frac{\delta}{a})\Bigl(\|\Gamma_k\rho\zeta \Psi^+\varphi)\| +
\|[\rho\zeta,\Gamma_k]\Psi^+\varphi)\|\Bigr)\\& \le 
& A_\delta
\exp(\frac{\delta}{a})\Bigl(2^{-k\varepsilon}\|\zeta\Gamma_k(\rho\zeta \Psi^+
\varphi)\|_\varepsilon
+2^{-k}C\|\Psi^+\varphi)\|\Bigr) \\& \le  &
A'_\delta 2^{(\frac{\delta}{a}-k\varepsilon)}\sqrt{Q_b(\zeta\Gamma_k\Psi^+\varphi,
\zeta\Gamma_k\Psi^+\varphi)} + 
A''_\delta 2^{(\frac{\delta}{a}-k)}\|\Psi^+\varphi)\|.
\end{eqnarray*}
\noindent
Set $a=\frac{2\delta}{\varepsilon k}$ with $k>K_\delta$, where $K_\delta$ is chosen large
enough so that 
\begin{eqnarray*}
A'_\delta 2^{(\frac{\delta}{a}-k\varepsilon)}&\le&\frac{\delta}{k}
\\  \noalign{\vskip-2pt}
\noalign{\noindent and}
 \noalign{\vskip-2pt}
A''_\delta 2^{(\frac{\delta}{a}-k)}&\le& \bigg(\frac{2}{3}\bigg)^k.
\end{eqnarray*}
Note that $\|\Psi^+\Gamma_k\varphi\|_{-1}\sim 2^{-k}\|\Psi^+\Gamma_k\varphi\|$. Thus from Lemma 4.2 we obtain
\noindent
\begin{eqnarray*} 
\|\Gamma_k\zeta \Psi^+\varphi\|^2 &= &\|\Gamma_k\zeta\Psi^+\varphi\|_{S_a}^2
+\|\Gamma_k\zeta \Psi^+\varphi\|_{M-S_a}^2\\ &\le &
C\Bigl(\frac{\delta^2}{k^2}Q_b(\zeta\Gamma_k\Psi^+\varphi,\zeta\Gamma_k\Psi^+\varphi)+
\Bigl(\frac{2}{3}\Bigr)^k\|\Psi^+\varphi)\|^2\Bigr)
\\&
\le &C\Bigl(\frac{\delta^2}{k^2}(\|\Gamma_k \bar{\partial}_b \zeta\Psi^+\varphi\|^2 +
\| \Gamma_k\bar{\partial}^*_b\zeta\Psi^+ \varphi\|^2)
+\Bigl(\frac{2}{3}\Bigr)^k\|\Psi^+\varphi\|^2\Bigr)
\end{eqnarray*}
\noindent
for $k>K_\delta$. Then, since $k \sim \log |\xi|$ for $\xi \in {\rm supp}(\gamma_k)$
we obtain, after multiplying by $k^2$ and summing over $k$,

\begin{eqnarray*}
\|(\log\Lambda)\zeta\Psi^+\varphi\|^2& \sim &\sum_{k\le
K_\delta}\|(\log\Lambda)\Gamma_k\zeta\Psi^+\varphi\|^2 + 
\sum_{k>K_\delta}k^2 \|\Gamma_k\zeta\Psi^+\varphi\|^2\\& \le 
&C\Bigl(K_{\delta}^2\|\Psi^+\varphi\|^2 + \delta^2(\|\Gamma_k\bar{\partial}_b
\zeta\Psi^+\varphi\|^2 +\| \Gamma_k\bar{\partial}^*_b
\zeta\Psi^+\varphi\|^2)\Bigr)\\&\le
&\delta^2Q_b(\zeta\Psi^+\varphi,\zeta\Psi^+\varphi) +  C_\delta
\|\Psi^+\varphi\|^2.
\end{eqnarray*}
This concludes the proof for the superlogarithmic estimate in the $+$ microlocalization;
the proof in the $-$ case is entirely analogous. \hfill\qed\vglue12pt

Combining the $+$, $-$, and $0$ microlocalizations we obtain the following.

\proclaim{{C}orollary} Let $M$ be a pseudoconvex manifold. Let $S\subset M$ be a 
manifold such that the holomorphic dimension of $S$ at each point is less than or equal to
$m${\rm ,} where $m=\min\{q-1,n-q-2\}$. Let $x_o\in S${\rm ,} let $U$ be a neighborhood of $x_o$ and
suppose that
$\rho\in C^\infty(U)$ with $\rho \in {\cal I}_q(x_o)$ satisfying $(*)$.  
Then the superlogarithmic estimate $({\rm SL})_q$ holds. For each $\delta > 0$ there exists
$C_\delta$ such that
$$
({\rm SL})_q\quad\quad\|(\log\Lambda)\varphi\|^2\le \delta^2Q_b(\varphi,\varphi) + 
C_\delta \|\varphi\|^2 ,
$$
for all  $\varphi\in {\cal A}^{0,q}$ with support in $U$. 
\endproclaim 

\vglue-12pt
\section{Hypoellipticity}

In this section we show that the superlogarithmic estimate $({\rm SL})_q$ implies 
hypoellipticity of $\Box_b$. Further we show that if the range of $\bar{\partial}_b$ is
closed in $L^2$ then the estimates $({\rm SL}^+)_1$ on $(p,1)$-forms and $({\rm SL}^-)_0$ on
$(p,0)$-forms imply that the restrictions of $\Box_b$ to $(p,0)$-forms
orthogonal to ${\cal H}^{p,0}$ and to $(p,n-1)$-forms orthogonal to
${\cal H}^{p,n-1}$ are hypoelliptic.

\proclaim{Theorem} Assume that $({\rm SL})_q$ holds in a neighborhood $U$ of $x_o\in
M$. Then if $\varphi$ is a square integrable $(p,q)$\/{\rm -}\/form   such that $\Box_b\varphi
=\alpha$ with $\alpha$ square integrable whose restriction to $U$ is in $C^\infty(U)$
then the restriction of $\varphi$ to $U$ is also in $C^\infty(U)$.
\endproclaim   

\demo{Proof} More precisely, we will show that, for any $\zeta_0,\zeta_1\in
 C^\infty_0(U)$ with $\zeta_1=1$ in a neighborhood of the support of $\zeta_0$, if
$\zeta_1\alpha\in H^s$ then $\zeta_0\varphi\in H^s$. To do this we   first prove the
following {\it a~priori} estimate: given $s$ there exists $C_s>0$ such that
$$
*\quad\quad\|\zeta_0\varphi\|_s\le C_s(\|\zeta_1\alpha\|_s+\|\varphi\|),
$$
for all $\varphi\in{\cal A}^{p,q}_b$. Let $\sigma\in C^\infty_0(U)$ such that
$\sigma=1$ in a neighborhood of ${\rm supp}(\zeta_0)$ and  $\zeta_1=1$ in a neighborhood of
${\rm supp}(\sigma)$. Let $\zeta'\in C^\infty_0(U')$ with $\zeta'=1$ in a neighborhood of
$\bar{U}$. We define the operator
$R^s$ by
$$
R^su(x)=\int e^{ix\cdot\xi}(1+|\xi|^2)^{\frac{s\sigma(x)}{2}}\hat{u}(\xi)d\xi,
$$
for $u\in C^\infty_0(U')$. Since the symbol of $(\Lambda^s-R^s)\zeta_0$ is zero,  
$$
\|\zeta_0\varphi\|_s\le\|R^s(\zeta_0\varphi)\|+
C\|\varphi\|\le\|[R^s,\zeta_0](\zeta_1\varphi)\|+\|R^s(\zeta_1\varphi)\|+C\|\varphi\|.
$$

From the calculus of pseudodifferential operators we conclude that 
$$
\|[R^s,\zeta_0](\zeta_1\varphi)\|\le C(\|R^{s-1}(\zeta_1\varphi)\|+\|\varphi\|)
$$
and
\begin{eqnarray*}
\|R^s(\zeta_1\varphi)\|&= &\|R^s(\zeta'\zeta_1\varphi)\|\\
&\le& \|\zeta'R^s(\zeta_1\varphi)\|+\|[R^s,\zeta'](\zeta_1\varphi)\|\\
&\le & \|\zeta'R^s(\zeta_1\varphi)\|+O(\|R^{s-1}(\zeta_1\varphi)\|+\|\varphi\|)\\
&\le & \|\zeta'R^s(\zeta_1\varphi)\|+C\|\varphi\|
\end{eqnarray*}
so that
$$
\|\zeta_0\varphi\|_s\le C(\|\log(\Lambda)\zeta' R^s(\zeta_1\varphi)\|+\|\varphi\|).
$$
Next, from $({\rm SL})_q$ with $\varphi$ replaced by $\zeta' R^s(\zeta_1\varphi)$, we obtain
$$
**\quad\|(\log\Lambda)R^s(\zeta_1\varphi)\|^2\le
\delta^2Q_b(\zeta' R^s(\zeta_1\varphi),\zeta' R^s(\zeta_1\varphi)) + 
C_\delta \|\varphi\|^2.
$$
The equation $\Box_b\varphi=\alpha$ is equivalent to
$$
Q_b(\varphi,\psi)=(\alpha,\psi),
$$
for all $\psi\in{\cal A}^{p,q}_b$. So we have
\begin{eqnarray*}
Q_b(\zeta' R^s(\zeta_1\varphi),\zeta' R^s(\zeta_1\varphi)) &=& 
Q_b(\varphi,\zeta_1(R^s)^*(\zeta')^2 R^s(\zeta_1\varphi))+ {\it error}\\ &=&
(\alpha,\zeta_1(R^s)^*(\zeta')^2 R^s(\zeta_1\varphi))+{\it error}\\ &=&
(\zeta' R^s(\zeta_1\alpha),\zeta' R^s(\zeta_1\varphi))+{\it error}.
\end{eqnarray*}
The ``error'' is given by $error=I+II$, where 
$$
I=([\bar{\partial}_b,\zeta'R^s\zeta_1]\varphi,\bar{\partial}_b\zeta'
R^s\zeta_1\varphi)+
([\bar{\partial}^*_b,\zeta'R^s\zeta_1]\varphi,\bar{\partial}^*_b\zeta'
R^s\zeta_1\varphi)
$$
and
\begin{eqnarray*}
II&=&(\bar{\partial}_b\varphi,[\zeta_1{R^s}^*\zeta',\ \bar{\partial}_b]\zeta'
R^s\zeta_1\varphi)+
(\bar{\partial}^*_b\varphi,\ [\zeta_1{R^s}^*\zeta',\ \bar{\partial}^*_b]\zeta'
R^s\zeta_1\varphi)\\ &=&
([[\zeta_1{R^s}^*\zeta',\ \bar{\partial}_b]^*,\ \bar{\partial}_b]\varphi,\ \zeta'
R^s\zeta_1\varphi)+
([[\zeta_1{R^s}^*\zeta',\ \bar{\partial}^*_b]^*,\ \bar{\partial}^*_b]\varphi,\ \zeta'
R^s\zeta_1\varphi)\\&&+\
([\zeta_1{R^s}^*\zeta',\ \bar{\partial}_b]^*\varphi,\ \bar{\partial}_b^*\zeta'
R^s\zeta_1\varphi)+
([\zeta_1{R^s}^*\zeta',\ \bar{\partial}^*_b]^*\varphi,\ \bar{\partial}_b\zeta'
R^s\zeta_1\varphi).
\end{eqnarray*}
By the Jacobi identity,  
$$
[\bar{\partial}_b,\zeta'R^s\zeta_1]=[\bar{\partial}_b,\zeta']R^s\zeta_1+
\zeta'[\bar{\partial}_b, R^s]\zeta_1+\zeta'R^s[\bar{\partial}_b, \zeta_1].
$$
Since the supports of the derivatives of $\zeta_1$ and of $\zeta'$ are disjoint from the
support of $\sigma$ the operator
$[\bar{\partial}_b,\zeta']R^s\zeta_1+\zeta'R^s[\bar{\partial}_b, \zeta_1]$ is bounded
in $L_2$. The principal symbol of $[\bar{\partial}_b, R^s]$ is bounded by 
$C(\log (1+|\xi|^2))(1+|\xi|^2)^\frac{s\sigma (x)}{2}$; hence
$$
\|\zeta'[\bar{\partial}_b, R^s]\zeta_1\varphi\| \le 
C(\|(\log \Lambda) R^s\zeta_1\varphi\|+\|\varphi\|).
$$
Arguing similarly we can bound all the terms in $I$ and $II$ and obtain
$$
error^2\le C(Q_b(\zeta' R^s(\zeta_1\varphi),\zeta' R^s(\zeta_1\varphi))+ 
\|(\log \Lambda )R^s\zeta_1\varphi\|^2+\|\varphi\|^2).
$$
After combining this with $**$ we see that ${\it error}^2$ gets multiplied by $\delta^2$ and thus, for
suitably small $\delta$, we obtain the {\it a~priori} estimate $*$. To conclude the proof of
the theorem we must pass from the {\it a~priori} estimate to showing regularity of the
solution. This can be done by applying a standard smoothing operator to $\zeta_1\varphi$
or by using the method of elliptic regularization as in [KN]. \enddemo

Microlocally we obtain the following.

\proclaim{Lemma} With the same notation as above{\rm ,} estimate $({\rm SL}^+)_q$ implies\/{\rm :}
$$
\|\Psi^+R^s\zeta_1\varphi\|\le
C(\|\Psi^+R^s\zeta_1\Box_b\varphi\|+\|\varphi\|)
$$
and
\begin{eqnarray*}
\|\Psi^+R^s\zeta_1\bar{\partial}_b\varphi\|+
\|\Psi^+R^s\zeta_1\bar{\partial}_b^*\varphi\|&\le&
C(\|\Psi^+R^s\zeta_1\Box_b\bar{\partial}_b\varphi\|+
\|\Psi^+R^s\zeta_1\Box_b\bar{\partial}_b^*\varphi\|\\
& &+\ \|\varphi\|),
\end{eqnarray*}
for all $\varphi\in{\cal A}^{p,q}_b\cap L_2$. Analogously{\rm ,} $({\rm SL}^-)_q$ implies\/{\rm :}\/
$$
\|\Psi^-R^s\zeta_1\varphi\|\le
C(\|\Psi^-R^s\zeta_1\Box_b\varphi\|+\|\varphi\|)
$$
and
\begin{eqnarray*}
\|\Psi^-R^s\zeta_1\bar{\partial}_b\varphi\|+
\|\Psi^-R^s\zeta_1\bar{\partial}_b^*\varphi\|&\le&
C(\|\Psi^-R^s\zeta_1\Box_b\bar{\partial}_b\varphi\|+
\|\Psi^-R^s\zeta_1\Box_b\bar{\partial}_b^*\varphi\|\\
& &+\ \|\varphi\|),
\end{eqnarray*}
for all $\varphi\in{\cal A}^{p,q}_b\cap L_2$.
\endproclaim

\demo{Proof} We have
\begin{eqnarray*}
Q_b(\zeta_1\Psi^+R^s\zeta'\varphi, \zeta_1\Psi^+R^s\zeta'\varphi) &= &
\|\Psi^+R^s\zeta_1\bar{\partial}_b\varphi\|^2+
\|\Psi^+R^s\zeta_1\bar{\partial}_b^*\varphi\|^2 + {\it error}\\ &= &
(\Psi^+\zeta_1\Box_b\varphi, \zeta_1\Psi^+\zeta'\varphi)+{\it error}\\ &\le &
C(\|\Psi^+\zeta_1\Box_b\varphi\|^2+\| \zeta_1\Psi^+\zeta'\varphi)\|^2)+{\it error}\\ &\le &
C(\|\Psi^+\zeta_1\bar{\partial}_b\bar{\partial}_b^*\varphi\|^2+
\|\Psi^+\zeta_1\bar{\partial}^*_b\bar{\partial}_b\varphi\|^2)+{\it error},
\end{eqnarray*}
and
\begin{eqnarray*}
\|\Psi^+\zeta_1\bar{\partial}^*_b\bar{\partial}_b\varphi\|^2
&=& (\Psi^+\zeta_1\bar{\partial}_b\bar{\partial}^*_b\bar{\partial}_b\varphi,
\Psi^+\zeta_1\bar{\partial}_b\varphi)+{\it error}\\&\le &
(\|\Psi^+\zeta_1\Box_b\bar{\partial}_b\varphi\|\|\Psi^+\zeta_1\bar{\partial}_b\varphi\|)
+{\it error},
\end{eqnarray*}
\begin{eqnarray*}
\|\Psi^+\zeta_1\bar{\partial}_b\bar{\partial}^*_b\varphi\|^2
&=& (\Psi^+\zeta_1\bar{\partial}_b\bar{\partial}^*_b\bar{\partial}_b\varphi,
\Psi^+\zeta_1\bar{\partial}^*_b\varphi)+{\it error}\\&\le &
(\|\Psi^+\zeta_1\Box_b\bar{\partial}^*_b\varphi\|\|\Psi^+\zeta_1\bar{\partial}^*_b\varphi\|)
+{\it error}.
\end{eqnarray*}
The ${\it error}$ terms arise from commutators such as $[\zeta_1\Psi^+R^s,\bar{\partial}_b]$.
These are ana\-lyzed as before except for the terms that involve the commutators
$[\Psi^+,\bar{\partial}_b]$ and $[\Psi^+,\bar{\partial}^*_b]$. To bound such terms let
$\tilde{\psi}^0$ be a symbol which equals one in a neighborhood of the support of the
derivatives of $\psi^+$ and whose support is contained in a region of the form 
$$
\{\xi=(\xi',\xi_{2n-1})\in{\Bbb R}^{2n-1}\big{|}\ A<\xi_{2n-1}<a|\xi'|\}
$$
and let $\tilde{\psi}^+$ be a symbol which equals one on the support of
$\tilde{\psi}^0$ with support contained in $\{\xi_{n-1}>0\}$. Denoting by
$\tilde{\Psi}^0$ and by $\tilde{\Psi}^+$ the corresponding pseudodifferential operators
we have
$$   
\|\tilde{\Psi}^0R^s\zeta_1\varphi\| \le
C(\|\Box_b\tilde{\Psi}^0R^s\zeta_1
\varphi\|_{-2}+\|\tilde{\Psi}^0R^s\zeta_1\varphi\|_{-2})
$$
since the support of $\tilde{\psi}^0$ lies in the elliptic region. These are
the terms that arise in the ${\it error}$. Thus we have
\begin{eqnarray*}
\|\tilde{\Psi}^0R^s\zeta_1\varphi\| &\le &
C(\|\tilde{\Psi}^+R^{s-2}\zeta_1\Box_b
\varphi\|+\|\tilde{\Psi}^+R^{s-2}\zeta_1\varphi\|+\|\varphi\|)\\&\le &
C(\|\tilde{\Psi}^+R^{s-2}\zeta_1\bar{\partial}_b\varphi\|+
\|\tilde{\Psi}^+R^{s-2}\zeta_1\bar{\partial}_b^*\varphi\|\\
&&+\
\|\tilde{\Psi}^+R^{s-2}\zeta_1\varphi\|+\|\varphi\|).
\end{eqnarray*}
We feed this into the above estimates with $\Psi^+$ replaced by
$\tilde{\Psi}^+$ and thus obtain the desired estimates in the $+$ microlocalization and
a parallel argument yields them in the $-$ microlocalization.\enddemo

The following result deals with the cases of $(p,0)$-forms and of $(p,n-1)$-forms,
where the spaces ${\cal H}^{p,0}_b$ and ${\cal H}^{p,n-1}_b$ are infinite-dimensional. We assume that $({\rm SL}^+)_1$
and
$({\rm SL}^-)_0$ hold. This assumption is equivalent to assuming that $({\rm SL}^-)_{n-2}$ and $({\rm SL}^+)_{n-1}$ hold,
as is seen from the following observation.

\demo{{R}emark} When we use  the operator $F_q$, as in the proof of Proposition 5.3, it follows that the estimate $({\rm
SL}^+)_q$ holds if and only if $({\rm SL}^-)_{n-q-1}$ holds.
\enddemo

\proclaim{Theorem} Let $M$ be a {\rm CR} manifold and assume that the operator $\bar{\partial}_b$ on
functions has closed range in $L_2(M)$. Assume that $({\rm SL}^+)_1$ and $({\rm SL}^-)_0$ hold on a
neighborhood $U$. Suppose that $q$ equals either $0$ or $n-1${\rm ,} 
$\varphi\perp{\cal H}^{p,q}_b${\rm ,} and $\Box_b\varphi=\alpha$ with $\alpha$ restricted
to $U$   in $C^\infty (U)$. Then the restriction of $\varphi$ is also in $C^\infty
(U)$.
\endproclaim

\demo{Proof} We will deal only with the case of functions, that is $(0,0)$-forms;
the remaining cases then follow. We will show that if $u\in L_2(M)$,
$u\perp{\cal H}^{0,0}_b$ and if
$\Box_bu=f$ with
$f\in L_2(M)$ and the restriction of $f$ to $U$ is in $C^\infty(U)$ then 
the restriction of $u$ to $U$ is in $C^\infty(U)$. Assuming that $u\in C^\infty(M)$ and
using the above lemma we obtain the following {\it a~priori} estimate: 
$$
\|\Psi^-R^s\zeta_1u\|\le
C(\|\Psi^-R^s\zeta_1f\|+\|u\|).
$$
Note that we have not used the assumption $u\perp{\cal H}^{0,0}_b$. We can now
deduce that $\|\Psi^-R^s\zeta_1u\|<\infty$ whenever $\Box_bu=f$. Since the range of
$\bar{\partial}_b$ is closed we conclude that the range of $\bar{\partial}^*_b$ is
closed and since on functions we have $\Box_bu=\bar{\partial}_b^*\bar{\partial}_bu$ 
we conclude that the range of $\bar{\partial}^*_b$ equals the
orthogonal complement of ${\cal H}^{0,0}_b$. Thus there exists a $(0,1)$-form
$\varphi$ such that $\bar{\partial}_b^*\varphi=u$ and we choose $\varphi$ so that it
also satisfies $\bar{\partial}_b\varphi=0$.  Hence 
$$
\|\Psi^+R^s\zeta_1\bar{\partial}_b^*\varphi\|\le
C(\|\Psi^+R^s\zeta_1\Box_b\bar{\partial}_b^*\varphi\|+ \|\varphi\|).
$$
By elliptic regularization we conclude that whenever the right-hand side above is
finite then so is the left-hand side. Since the range of $\bar{\partial}^*_b$ is
closed,
$\|\varphi\|\le C\|\bar{\partial}_b^*\varphi\|=C\|u\|$. Hence
$$
\|\Psi^+R^s\zeta_1u\|\le C(\|\Psi^+R^s\zeta_1f\|+ \|u\|).
$$
Again we conclude that $\|\Psi^+R^s\zeta_1u\|<\infty$ and combining this with the above
and with the $0$ microlocalization we have $\zeta_0u\in H^s$ for all $s$, which
concludes the proof. \enddemo

{\it Remark}. The range of $\bar{\partial}_b$ on functions is closed whenever $M$
is compact and $({\rm SL})_1$ holds at all points of $M$. As noted this condition cannot
hold if ${\rm dim}(M)=3$. Another condition for the closed range of $\bar{\partial}_b$ is
given in [K4], namely: if $X$ is an $n$-dimensional Stein manifold and if 
$\Omega\subset X$ is a relatively compact pseudoconvex domain with a smooth boundary $M$ then 
$\bar{\partial}_b$ on $M$ has closed range (see [K4]).
 
\section{The $\bar{\partial}$-Neumann problem}

In this section we shall establish the relation between subellipticity, subelliptic multipliers, and superlogarithmic
estimates between boundaries of pseudoconvex domains and the corresponding estimates for domains. These results can be
extended to various other estimates. The method here is the type of microlocalization worked out in [K4]. Another
reduction to the boundary was derived by Greiner and Stein (see [GS]).

Let $\Omega\subset X$ be a relatively compact domain with a smooth boundary in a complex hermitian manifold
$X$ and let $M$ denote the boundary of $\Omega$. Let $r$ be a defining function for $M$; that is, $r$ is a $C^\infty$
function in a neighborhood of $M$ such that $r=0$ on $M$ and $dr\neq 0$. We will further assume that $|dr|=1$ on
$M$, that $r<0$ in $\Omega$, and that $r>0$ outside of $\bar{\Omega}$. To simplify
formulas we will also assume that $|r|$ is the geodesic distance to M. We denote by
${\cal A}^{p,q}$ the space of $(p,q)$-forms in $C^\infty(X)$ restricted to $\bar{\Omega}$. As usual,
$\bar{\partial}:{\cal A}^{p,q}\to{\cal A}^{p,q+1}$. We denote by
$\bar{\partial}^*$ the $L_2$-adjoint of $\bar{\partial}$. Also, $\bar{\partial}^*$ is an 
unbounded operator in Hilbert space whose domain, denoted by ${\rm Dom}(\bar\partial^*)$, consists of all $\varphi\in
L_2^{p,q}(M)$ such that $(\varphi,\bar{\partial}u)=(\bar{\partial}^*\varphi,u)$ for all
$u\in L_2^{p,q-1}(M)$ with $\bar{\partial}u\in L_2^{p,q}(M)$. If $x_0\in M$ we  choose 
local holomorphic coordinates $\{z_1,z_2,\ldots ,z_n\}$ on a neighborhood $U$ of $x_0$ with
origin at $x_0$ such that $dz_n|_0=\partial r|_0$. Let $\{\omega_1,\ldots ,\omega_n\}$ be an
orthonormal basis of the $(1,0)$-forms on $U$ with  $\omega_n=\partial r$, so that $L_n(r)=1$. Let $\{L_1,\ldots ,L_n\}$ be the
dual basis to the $\omega_i$ and note that on $M$ we have $L_ir=0$ for $i=1,\ldots ,n-1$, so that the CR structure
on $M$, and on the manifolds $r=constant$, is given by the $\{L_1,\ldots ,L_{n-1}\}$.
Setting $T=\frac{1}{\sqrt 2}(L_n-\bar{L}_n)$, we have that
$\{L_1,..L_{n-1},\bar{L}_1,\ldots ,\bar{L}_{n-1},T\}$ is an orthonormal
basis of vectors tangent to $M$. Also, $T$ satisfies $T(r)=0$ and $\bar{T}=-T$.  

Suppose that $\varphi\in{\cal A}^{0,q}\cap C_0^\infty(U\cap\bar{\Omega})$, so that
$\varphi=\sum\varphi_I\bar{\omega}_I$. Then $\varphi$ is in the domain of
$\bar{\partial}^*$ if and only if $\varphi_I=0$ on $M$ whenever $n\in I$. The
$\bar{\partial}$-Neumann problem consists of solving the equation $\Box\varphi=\alpha$,
where $\Box= \bar{\partial}\bar{\partial}^*+\bar{\partial}^*\bar{\partial}$. Setting
$Q(\varphi,\psi)=(\bar{\partial}\varphi,\bar{\partial}\psi)+
(\bar{\partial}^*\varphi,\bar{\partial}^*\psi)$ we see that the above equation is
satisfied if and only if $\varphi$ is in the domain of $\bar{\partial}^*$ and
$Q(\varphi,\psi)=(\alpha,\psi)$ for all $\psi$ in the domain of $\bar{\partial}^*$.

On $U\cap\bar{\Omega}$ we will use the coordinates
$\{x_1,\ldots ,x_{2n-1},r\}$, where 
$x_i={\rm Re}(z_i)$ for $i=1,\ldots n-1$ and $x_i=\Im(z_{i-n+1})$ for $i=n,\dots,2n-1$ . We will denote by
$x=(x_1,\ldots ,x_{2n-1})$ the tangential coordinates and by
$\xi=(\xi_1,\ldots ,\xi_{2n-1})$ the dual coordinates. We define the tangential Fourier
transform of $u\in\ C^\infty_0(U\cap\bar{\Omega}$ by
$$
\tilde{u}(\xi,r)=\int e^{ix\cdot\xi}u(x,r)dx.
$$
This notation should not be confused with the one used in previous sections.
We define the operator $\Lambda^s_{\tan}$ by
$$
\widetilde{\Lambda^s_{\tan}u}(\xi,r)=(|\xi|^2+1)^{s/2}\tilde{u}(\xi,r)
$$
and the tangential Sobolev norms by
$$
|\!|\!|u|\!|\!|^2_s=\int |\Lambda^s_{\tan}u|^2d\xi dr .
$$
In terms of these coordinates the operators $L_i$ can be written as
$$
L_i=\sum a^k_i(x,r)\frac{\partial}{\partial x_k},
$$
for $i=1,\ldots ,n-1$ and
$$
L_n=\frac{\partial}{\partial r}+\sum a^k(x,r)\frac{\partial}{\partial x_k}.
$$ 
We define the tangential symbols $\mu_i$ by
$$
\mu_i(x,\xi,r)=\frac{1}{\sqrt{-1}}\sum a^k_i(x,r)\xi_k,
$$
for $i=1,\ldots ,n-1$ and
$$
\mu_n(x,\xi,r)=\frac{1}{\sqrt{-1}}\sum (a^k(x,r)-\bar{a}^k(x,r))\xi_k.
$$
Note that $\mu_n$ is real. We set $\dot{\mu}_i(x,\xi)=\mu_i(x,\xi,0)$ and
$\dot{\mu}=(\dot{\mu}_1,\ldots ,\dot{\mu}_n)$ for $i=1,\ldots ,n$, and
$$
|\dot{\mu}(x,\xi)|=\sqrt{\sum|\dot{\mu}_i(x,\xi)|^2} .
$$
Subelliptic multipliers for the $\bar{\partial}$-Neumann problem are defined as follows.

\numbereddemo{Definition} Suppose that $x_0\in M$ and that $f\in C^\infty(U)$.  Then $f$ is a
subelliptic multiplier on $(p,q)$-forms for the $\bar{\partial}$-Neumann problem
if and only if there exist $\varepsilon > 0$ and $C > 0$ such that
$$
|\!|\!|f\varphi|\!|\!|_\varepsilon^2\le C(Q(\varphi,\varphi)+\|\varphi\|^2),
$$   
for all ${\rm Dom}(\bar\partial^*)\cap\varphi\in{\cal A}^{p,q}\cap C^\infty_0(U\cap\bar{\Omega})$. We denote by
${\cal J}_q(x_0)$ the {\it ideal of germs of subelliptic multipliers} for the
$\bar{\partial}$-Neumann problem at $x_0$. 
\enddemo

The relation between the ideals of subelliptic multipliers for CR manifolds and for
the $\bar{\partial}$-Neumann problem is given in the following.
 
\proclaim{Theorem} Suppose that $\Omega\subset X$ is a pseudoconvex domain in a hermitian
manifold $X$. Suppose that $\Omega$ is compact and has a smooth boundary~$M$.
Then{\rm ,} if $x_0\in M${\rm ,} the ideals ${\cal J}_q(x_0)$ and ${\cal I}^+_q(x_0)$ are
related as follows. If $f\in{\cal J}_q(x_0)$ then the restriction of $f$ to $M$ is
in ${\cal I}^+_q(x_0)$. Furthermore{\rm ,} if $\rho\in{\cal I}^+_q(x_0)${\rm ,} if $f\in
C^\infty(\bar{\Omega})${\rm ,} and if the restriction of $f$ to $M$ equals $\rho$ then
$f\in{\cal J}_q(x_0)$.
\endproclaim

The fact that hypoellipticity of the $\bar{\partial}$-Neumann problem follows from the
assumption that a superlogarithmic estimate holds is formulated as follows.

\proclaim{Theorem} As above with $x_0\in M$ and $M$ the boundary of $\Omega,$ assume that
condition $({\rm SL}^+)_q$ holds in a neighborhood $U\cap M$ of $x_0$. Then if $\alpha$ is a
$(p,q)$\/{\rm -}\/form on $\Omega$ which is in  $L_2${\rm ,} whose the restriction to
$U\cap\bar{\Omega}$ is in $C^\infty(U\cap\bar{\Omega})$ and if $\varphi$ satisfies the
equation $\Box\varphi=\alpha$ {\rm (}\/this means{\rm ,} in particular{\rm ,} that $\varphi$ and
$\bar{\partial}\varphi$ are in the domains of $\bar{\partial}^*$ on $(p,q)$ and
$(p,q+1)$\/{\rm -}\/forms{\rm ,} respectively{\rm ),} then the restriction of $\varphi$ to
$U\cap\bar{\Omega}$ is in $C^\infty(U\cap\bar{\Omega})${\rm . }
\endproclaim  

The key to proving these theorems is a passage between forms in ${\cal A}^{p,q}_b$
and forms in ${\cal A}^{p,q}$. This is done by introducing an ``approximate'' harmonic
extension of $u$ on $U\cap M$ to $U\cap\bar\Omega$, denoted by $u^{(h)}$. Supposing that $u\in C^\infty_0(U\cap M)$ we
define 
$u^{(h)}\in C^\infty(\{(x,r)\in{\Bbb R}^{2n} |\  r\le 0\})$ by
$$
 u^{(h)}(x,r)=(2\pi)^{-2n+1}\int e^{ix\cdot\xi}e^{r|\dot{\mu}(x,\xi)|}\tilde{u}(\xi)d\xi,
$$
so that $u^{(h)}(x,0)=u(x)$.

In this section $\|\ \|$ will denote the $L_2$ norm on $\Omega$, $\|\ \|^b$ the $L_2$ norm on $M$,
and $\|\ \|^b_s$ the Sobolev $s$-norm on $M$.

\proclaim{Lemma} For each $k\in{\Bbb Z}${\rm ,} $k\ge 0${\rm ,} and $s\in{\Bbb R}$ there exists $C_{s,k}>0$
such that
$$
|\!|\!|r^ku^{(h)}|\!|\!|_s\le C_{s,k}\|u\|^b_{s-k-\frac{1}{2}},
$$
for all $u\in C^\infty_0(U\cap\bar{\Omega})$.
\endproclaim

\demo{Proof} We have
$$ 
|\!|\!|r^ku^{(h)}|\!|\!|_s^2\le C\int
r^{2k}e^{2r|\dot{\mu}(x,\xi)|}(1+|\xi|^2)^s|\tilde{u}(\xi)|^2d\xi dr.
$$
Substituting $r'= r|\dot{\mu}(x,\xi)|$ we integrate first with respect to $r'$, and
note that $|\dot{\mu}(x,\xi)|^2\sim (1+|\xi|^2)$ to obtain the result.\enddemo

Setting $\ \triangle=-\sum\frac{\partial^2}{\partial{z_i}\partial\bar{z}_i}\ $ we have
\begin{eqnarray*}
\triangle &=&-\sum_{i=1}^{n}
L_i\bar{L}_i+\sum_{i=1}^{2n-1}a^i(x,r)\frac{\partial}{\partial
x_i}+a(x,r)\frac{\partial}{\partial r}\\&=&
-\frac{\partial^2}{\partial^2r}+T^2-\sum_{i=1}^{n-1} L_i\bar{L}_i
+\sum_{i=1}^{2n-1}b^i(x,r)\frac{\partial}{\partial x_i}+b(x,r)\frac{\partial}{\partial
r},
\end{eqnarray*}
since $\frac{\partial}{\partial r}=\frac{1}{2}(L_n+\bar{L}_n)$ and 
$L_n\bar{L}_n=\frac{\partial^2}{\partial^2r}-T^2+D,$ where $D$ is a first order operator. Hence if $(x,r)\in
U\cap\bar{\Omega}$,  
$$
\triangle(u^{(h)})(x,r)=\int
e^{ix\cdot\xi}e^{r|\dot{\mu}(x,\xi)|}(p^1(x,r,\xi)+rp^2(x,r,\xi))\tilde{u}(\xi)d\xi+Eu(x,r),
$$
where the $p^k(x,r,\xi)$ are symbols of order $k$, uniformly in $r$.  Note that $E(u)$ denotes the
error term which is an operator of order $-\infty$. Abusing notation we will denote all such
terms by $E(u)$. Further, 
$$
L_iu^{(h)}(x,r)=(L_iu)^{(h)}+K_iu(x,r)+Eu(x,r),
$$
with
$$
K_iu(x,r)=\int
e^{ix\cdot\xi}e^{r|\dot{\mu}(x,\xi)|}(\dot{\mu}_i(x,\xi)+p_i^0(x,r,\xi)+
rp_i^1(x,r,\xi))\tilde{u}(\xi)d\xi
$$
and
$$
\bar{L}_iu^{(h)}(x,r)=(\bar{L}_iu)^{(h)}+\bar{K}_iu(x,r)+Eu(x,r)
$$
 for $i=1,\ldots ,n-1$. Since $L_n=\frac{1}{\sqrt{2}}(\frac{\partial}{\partial r}+T)$  
\begin{eqnarray*}
L_nu^{(h)}(x,r)&=&\frac{1}{\sqrt{2}}\int
e^{ix\cdot\xi}e^{r|\dot{\mu}(x,\xi)|}\Bigl(|\dot{\mu}(x,\xi)|+\dot{\mu}_n(x,\xi)+
p^0_n(x,r,\xi)\\&\ &\ \ +\ rp^1_n(x,r,\xi)\Bigr)\tilde{u}(\xi)d\xi\ +Eu(x,r)
\end{eqnarray*}
and
\begin{eqnarray*}
\bar{L}_nu^{(h)}(x,r)&=&\frac{1}{\sqrt{2}}\int
e^{ix\cdot\xi}e^{r|\dot{\mu}(x,\xi)|}\Bigl(|\dot{\mu}(x,\xi)|-\dot{\mu}_n(x,\xi)+
\bar{p}^0_n(x,r,-\xi)\\&\ &\ \ +\ r\bar{p}^1_n(x,r,-\xi)\Bigr)\tilde{u}(\xi)d\xi\ +Eu(x,r)
\end{eqnarray*}where the $p^k_i$ are symbols of order $k$. If $v\in
C^\infty_0(U\cap\bar\Omega)$ we define $\dot{v}$ to be the restriction of $v$ to $M$ and
we have
$$
\|\dot{v}\|^b_s\le C\lrp{|\!|\!|v|\!|\!|_{s+\frac{1}{2}}+
|\!|\!|\frac{\partial v}{\partial r}|\!|\!|_{s-\frac{1}{2}}}.
$$

We microlocalize $v$ for each fixed $r$, setting
$$
\widetilde{\Psi v}(\xi,r)=\psi(\xi)\tilde{v}(\xi,r)
$$
and then we have:

\proclaim{Lemma} If $U$ is sufficiently small then there exists $C>0$ such that
$$
\|\Psi^+\bar{L}_n\dot{v}^{(h)}\|\le
C\lrp{\sum_{1}^{n-1}\|\Psi^+\bar{L}_iv\|^b_{-\frac{1}{2}}+\|v\|}
$$
and
$$
|\!|\!|\Psi^0\dot{v}^{(h)}|\!|\!|_1+|\!|\!|\Psi^-\dot{v}^{(h)}|\!|\!|_1\le
C\lrp{\sum_{1}^{n}\|\bar{L}_iv\|+\|v\|},
$$
for all $v\in C^\infty_0(U\cap\bar{\Omega})$.
\endproclaim

\demo{Proof} Choosing $U$ sufficiently small we have $\dot{\mu}_n(x,\xi)\ge 0$ when
$\xi\in {\rm supp}(\psi^+)$ and $\dot{\mu}_n(x,\xi)\le 0$ when $\xi\in {\rm supp}(\psi^-)$. Then
for $x\in U$ and $\xi\in {\rm supp}(\psi^+)$,  
$$
|\dot{\mu}(x,\xi)|-\dot{\mu}_n(x,\xi)=
\sum_{1}^{n-1}\frac{\dot{\mu}_i(x,\xi)}
{|\dot{\mu}(x,\xi)|+\dot{\mu}_n(x,\xi)}\dot{\overline{\mu}}_i(x,\xi);
$$
hence the right-hand side is the principal symbol of a pseudodifferential
operator on $M$ of the form $\sum_{1}^{n-1}P_i\bar{L}_i$, where the $P_i$ are of order
zero. This establishes the first inequality. The second inequality follows from the
fact that for $x\in U$ and $\xi\in {\rm supp}(\psi^-)$,
$$
|\dot{\mu}(x,\xi)|-\dot{\mu}_n(x,\xi)\ \ge\  C|\xi|,
$$
when $U$ is sufficiently small.\enddemo

If $\varphi\in{\cal A}^{p,q}$ we denote by  $\dot{\varphi}$ the restriction of
$\varphi$ to $M$. Note also that if  $\varphi$ is in the domain of $\bar{\partial}^*$
then $\dot{\varphi}\in{\cal A}^{p,q}_b$. Recall that on a pseudoncovex domain
the following are equivalent:
$$
Q(\varphi,\varphi)+\|\varphi\|^2\sim\sum\int_Mc_{IJ}\varphi_I\bar{\varphi}_JdS+
\sum_{1}^n\|\bar{L}_i\varphi\|^2+\|\varphi\|^2,
$$
for all $\varphi\in{\cal A}^{p,q}\cap C^\infty_0(U'\cap \bar\Omega)$ intersected with the
domain of $\bar{\partial}^*$. Also
\begin{eqnarray*}
&&\hskip-24pt Q_b(\zeta\Psi^+\varphi,\zeta\Psi^+\varphi)+(\|\zeta\Psi^+\varphi\|^b)^2\\
&&\quad \sim\
\sum(c_{IJ}\Lambda^\frac{1}{2}\zeta\Psi^+\varphi_I,
\Lambda^\frac{1}{2}\zeta\Psi^+\varphi_J)^b
 + \sum_{1}^{n-1}(\|\bar{L}_i\zeta\Psi^+\varphi\|^b)^2+(\|\zeta'\Psi'^+\varphi\|^b)^2,
\end{eqnarray*}
for all $\varphi\in{\cal A}_b^{p,q}$ with support in $U'\cap M$. Combining the
above with Lemma 8.5, we obtain:

\proclaim{Lemma} Suppose that $\Omega\subset X$ is a pseudoconvex domain in a hermitian
manifold $X$. Suppose that $\Omega$ is compact and has a smooth boundary $M$.
Then{\rm ,} if $x_0\in M$ and if $U$ is a sufficiently small neighborhood of $x_0$ then
there exists a constant $C>0$ such that
$$
|\!|\!|\Psi^-\varphi|\!|\!|_1^2+|\!|\!|\Psi^0\varphi|\!|\!|_1^2\le
C(Q(\varphi,\varphi)+\|\varphi\|^2),
$$
for all $\varphi\in{\cal A}^{p,q}\cap \{domain\  of\  \bar{\partial}^*\}$ with
support in
$U\cap\bar{\Omega}$.
\endproclaim

\demo{Proof} We have
$$
|\!|\!|\Psi^-\dot{\varphi}^{(h)}|\!|\!|_1^2+|\!|\!|\Psi^0\dot{\varphi}^{(h)}|\!|\!|_1^2\le
C\lrp{\sum_{1}^{n}\|\bar{L}_i\varphi_I\|^2+\|\varphi\|^2)\le 
C(Q(\varphi,\varphi)+\|\varphi\|^2}.
$$
Let $\varphi^{(0)}=\varphi-\dot{\varphi}^{(h)},\ $  so that
$\varphi^{(0)}=0$ on $M$. Then, with $\zeta\in C^\infty_0(U')$ and $\zeta=1$ on
a neighborhood of $\bar{U}$, 
\begin{eqnarray*}
|\!|\!|\Psi^-\varphi^{(0)}|\!|\!|_1^2+|\!|\!|\Psi^0\varphi^{(0)}|\!|\!|_1^2 &\le & 
C(Q(\zeta\Psi^-\varphi^{(0)},\zeta\Psi^-\varphi^{(0)})\\&\ &+
Q(\zeta\Psi^0\varphi^{(0)},\zeta\Psi^0\varphi^{(0)})+\|\varphi\|^2)\\&\le &
C(Q(\varphi,\varphi)+\|\varphi\|^2),
\end{eqnarray*}
which concludes the proof.\enddemo

Now we are in a position to prove the theorems of this section.

\demo{Proof of Theorem {\rm 8.2}} Suppose that $f\in{\cal J}_q(x_0)$ and
$\dot{f}=\rho$. Then, if $\varphi\in{\cal A}_b^{p,q}\cap
C^\infty_0(U'\cap\bar{\Omega})$,  
\begin{eqnarray*}
(\|\rho\zeta\Psi^+\varphi\|^b_\varepsilon )^2 &\le &
C(|\!|\!|f\zeta\Psi^+\varphi^{(h)}|\!|\!|_{\varepsilon+\frac{1}{2}}^2+
|\!|\!|f\zeta\Psi^+\frac{\partial\varphi^{(h)}}{\partial r}|\!|\!|_{\varepsilon-\frac{1}{2}}^2)
\\ &\le &
C(Q(\zeta\Lambda^{\frac{1}{2}}_{{\rm tan}}\Psi^+\varphi^{(h)},
\zeta\Lambda^{\frac{1}{2}}_{{\rm tan}}\Psi^+\varphi^{(h)})+
\|\Lambda^{\frac{1}{2}}_{{\rm tan}}\Psi^+\varphi^{(h)}\|^2)\\&\le &
C(Q_b(\zeta\Psi^+\varphi,\zeta\Psi^+\varphi)+(\|\Psi^+\varphi\|^b)^2),
\end{eqnarray*}
so that $\rho\in{\cal I}^+_q(x_0)$.

Next, assume that $\rho\in{\cal I}^+_q(x_0)$ and that $\dot{f}=\rho$. Then
if 
$$\varphi\in{\cal A}^{p,q}\cap
C^\infty_0(U\cap\bar{\Omega})\cap\{domain\ of\ \bar{\partial}^*\},
$$
 we have
\begin{eqnarray*}
|\!|\!|f\zeta\Psi^+\varphi|\!|\!|_\varepsilon^2 &\le & 
C(|\!|\!|f\zeta\Psi^+\dot{\varphi}^{(h)}|\!|\!|_\varepsilon^2+
|\!|\!|f\zeta\Psi^+\varphi^0|\!|\!|_\varepsilon)\\ &\le &
C(\|\rho\zeta\Psi^+\dot{\varphi}\|^b_{\varepsilon-\frac{1}{2}}+\|\varphi\|^2)\\&\le&
C(Q_b(\zeta\Lambda^{-\frac{1}{2}}\Psi^+\dot{\varphi},
\zeta\Lambda^{-\frac{1}{2}}\Psi^+\dot{\varphi})+\|\varphi\|^2)\\ &\le &
C(\sum_I\sum_{1}^{n-1}(\|\bar{L}_i\zeta\Lambda^{-\frac{1}{2}}
\Psi^+\dot{\varphi}\|^b)^2\\&\ &+\sum_{IJ}(c_{IJ}T\zeta\Lambda^{-\frac{1}{2}}
\Psi^+\dot{\varphi}_I,\zeta\Lambda^{-\frac{1}{2}}\Psi^+\dot{\varphi}_J)^b
+ \|\varphi\|^2).
\end{eqnarray*}
To estimate the last two terms above we proceed as follows. For $i<n$,
\begin{eqnarray*}
(\|\bar{L}_i\zeta\Lambda^{-\frac{1}{2}}
\Psi^+\dot{\varphi}\|^b)^2&=&2(\frac{\partial}{\partial
r}\bar{L}_i\zeta\Lambda^{-\frac{1}{2}}_{{\rm tan}}\varphi,
\bar{L}_i\zeta\Lambda^{-\frac{1}{2}}_{{\rm tan}}\varphi)\\&=&
2(\frac{\partial}{\partial r}\Lambda^{-1}_{{\rm tan}}\bar{L}_i\varphi,\bar{L}_i\varphi)+ {\it error}\\&\le &
C(|(\Lambda^{-1}_{{\rm tan}}\bar{L}_i\bar{L}_n\varphi,\bar{L}_i\varphi)|+
|(\Lambda^{-1}_{{\rm tan}}T\bar{L}_i\varphi,\bar{L}_i\varphi)|+{\it error}\\&\le &
C(Q(\varphi,\varphi)+\|\varphi\|^2),
\end{eqnarray*}
since $\Lambda^{-1}_{{\rm tan}}\bar{L}_i$ and $\Lambda^{-1}T$ are tangential pseudodifferential operators
of order zero. The error terms can be estimated since $[\frac{\partial}{\partial r},\bar{L}_i]$ is
a tangential first order operator. Finally, setting $P=(\zeta)^2T\Lambda^{-1}_{{\rm tan}}(\Psi^+)^2$ we
have
\begin{eqnarray*}
\sum_{IJ}(c_{IJ}T\zeta\Lambda^{-\frac{1}{2}}
\Psi^+\dot{\varphi}_I,\zeta\Lambda^{-\frac{1}{2}}\Psi^+\dot{\varphi}_J)^b &=&
\sum_{IJ}(c_{IJ}P\dot{\varphi}_I,\dot{\varphi}_J)^b+{\it error}\\ &\le &
|\sum_{IJ}(c_{IJ}P\dot{\varphi}_I,P\dot{\varphi}_J)^b|\\&\ &+
|\sum_{IJ}(c_{IJ}\dot{\varphi}_I,\dot{\varphi}_J)^b|+{\it error}\\ &\le &
C(Q(P\varphi,P\varphi)+Q(\varphi,\varphi)+\|\varphi\|^2)\\ &\le &
C(Q(\varphi,\varphi)+\|\varphi\|^2),
\end{eqnarray*} 
which concludes the proof of the theorem.\enddemo

\demo{Proof of Theorem {\rm 8.3}} Assume that the condition $({\rm SL}^+)_q$ holds in $U\cap
M$. Let $\varphi$ be a $(p,q)$-form that satisfies $\Box\varphi=\alpha$ with
$\alpha\in L_2$ and the restriction of $\alpha$ to $U\cap\bar{\Omega}$ in
$C^\infty(U\cap\bar{\Omega})$. We will show that the restriction of $\varphi$  to 
$U\cap\bar{\Omega}$ is in $C^\infty(U\cap\bar{\Omega})$. First we establish the
following {\it a~priori\/} estimate for the tangential derivatives of $\varphi$. Given
$s\in{\Bbb R}$ and $\zeta_0, \zeta_1\in C^\infty_0(U\cap\bar{\Omega})$  with
$\zeta_1=1$ on the support of $\zeta_0$, there exists a constant $C_s$ such that
$$
(*)_{{\rm tan}}\quad\quad|\!|\!|\zeta_0\varphi|\!|\!|_s\le C_s(|\!|\!|\zeta_1\alpha|\!|\!|_s+\|\varphi\|),
$$
for all $\varphi\in{\cal A}^{p,q}\cap\{domain\ of\ \bar{\partial}^*\}$. Now  let
$\zeta'\in C^\infty_0(U'\cap\bar{\Omega})$ with $\zeta'=1$ on $U\cap\Omega$.  We will
abuse notation and  use $\varphi$ to denote both $\zeta'\varphi$ and $\varphi$;
it will be clear from the context which is which and the errors committed will be
controlled as in Section 7. Again we set $\varphi=\dot{\varphi}^{(h)}+\varphi^{(0)}$. Since $\varphi^0=0$ on $M$,
$$ 
|\!|\!|\zeta_0\varphi^{(0)}|\!|\!|_{s+2}\le
C_s(|\!|\!|\zeta_1\Box\varphi^{(0)}|\!|\!|_s+\|\varphi^{(0)}\|)
$$
and
$$
\Box\varphi^{(0)}=\Box\varphi-\Box\dot{\varphi}^{(h)}=\alpha +P_1\varphi,
$$
where $P_1$ is a tangential pseudodifferential operator of order one. Since the above
inequality holds for all suitable pairs $\zeta_0, \zeta_1$ we deduce that
$$ 
|\!|\!|\zeta_0\varphi^{(0)}|\!|\!|_{s+2}\le C_s(|\!|\!|\zeta_1\alpha|\!|\!|_s+\|\varphi\|).
$$
Then, Lemma 8.6 implies that
$$
|\!|\!|\Psi^-\zeta_0\varphi|\!|\!|_{s+2}+|\!|\!|\Psi^0\zeta_0\varphi|\!|\!|_{s+2}\le
C_s(|\!|\!|\alpha|\!|\!|_s+\|\varphi\|).
$$
Hence, in order to prove $(*)_{\rm tan}$, it will suffice to prove
$$
(**)_{\rm tan}\quad\quad |\!|\!|\Psi^+\zeta_0\dot{\varphi}^{(h)}|\!|\!|_s\le 
C_s(|\!|\!|\alpha|\!|\!|_s+\|\varphi\|). 
$$
Now, following the argumentation of Section 7, in order to simplify the
formulas, we denote by $C$ and $C_\delta$ constants which may differ in different lines:
\begin{eqnarray*}
|\!|\!|\Psi^+\zeta_0\dot{\varphi}^{(h)}|\!|\!|^2_s&\le &
C\Bigl( (\|\Psi^+\zeta_0\dot{\varphi}\|_{s-\frac{1}{2}}^b)^2+\|\varphi\|^2\Bigr)\\&\le
& C\Bigl(
(\|\log(\Lambda_{\rm tan})R^{s-\frac{1}{2}}\Psi^+\dot{\varphi}\|^b)^2+\|\varphi\|^2\Bigr)
\\&\le & \delta^2Q_b(\zeta_1R^{s-\frac{1}{2}}\Psi^+\dot{\varphi},
\zeta_1R^{s-\frac{1}{2}}\Psi^+\dot{\varphi})+C_\delta\|\varphi\|^2\\&\le &
C\delta^2\Bigl(
\sum_{1}^{n-1}(\|\bar{L}_i\zeta_1R^{s-\frac{1}{2}}\Psi^+\dot{\varphi}\|^b)^2\\&  &+\
\sum_{IJ}(c_{IJ}T\zeta_1R^{s-\frac{1}{2}}\Psi^+\dot{\varphi}_I,
\zeta_1R^{s-\frac{1}{2}}\Psi^+\dot{\varphi}_J)^b\Bigr)+
C_\delta\|\varphi\|^2\\&\le &
C\delta^2\Bigl(
\sum_{1}^{n}\|\bar{L}_i\zeta_1R^s\Psi^+\varphi\|^2\\&  &+\
\sum_{IJ}\int_Mc_{IJ}T\zeta_1R^s\Psi^+\varphi_I
\zeta_1R^s\Psi^+\bar{\varphi}_JdS\Bigr)+
C_\delta\|\varphi\|^2\\&\le &
C\delta^2Q(\zeta_1R^s\Psi^+\varphi,\zeta_1R^s\Psi^+\varphi)+C_\delta\|\varphi\|^2
\\&\le & C\delta^2(R^s\Psi^+\zeta_1\alpha,\zeta_1R^s\Psi^+\varphi)+
C_\delta\|\varphi\|^2.
\end{eqnarray*}
The estimate $(*)_{\rm tan}$ then follows and  implies the next inequalities:
$$
|\!|\!|\frac{\partial}{\partial r}(\zeta_0\varphi)|\!|\!|_{s-1}\le 
C(|\!|\!|\zeta_1\alpha|\!|\!|_s+\|\varphi\|)
$$
and
$$
|\!|\!|\frac{\partial^{k+2}}{\partial r^{k+2}}(\zeta_0\varphi)|\!|\!|_{s-k-2}\le
C\lrp{\sum_{0}^{k}|\!|\!|\frac{\partial^j}{\partial r^j}(\zeta_1\alpha)|\!|\!|_{s-j}+\|\varphi\|},
$$
for integers $k\ge 0$. The first inequality follows from
\begin{eqnarray*}
\|\frac{\partial}{\partial r}(\zeta_0\varphi)\|^2&\le &\|\bar{L}_n(\zeta_0\varphi)\|^2
+\|T(\zeta_0\varphi)\|^2\\&\le &C(Q(\zeta_0\varphi,\zeta_0\varphi)+
|\!|\!|\zeta_0\varphi|\|_1^2+\|\varphi\|^2),
\end{eqnarray*}
and  is then obtained by substituting $\Lambda_{\rm tan}^{s-1}\varphi$ for
$\varphi$. The second inequality is obtained as follows. Since $\Box$ is elliptic we can solve the equation
$\Box\varphi=\alpha$ for the second derivatives with respect to $r$:
$$
\frac{\partial^2\varphi_I}{\partial r^2}=\sum_Ka^K_I\alpha_K+
\sum_{K,i,j}b^{Kij}_I\frac{\partial^2\varphi_K}{\partial x_i\partial x_j}+
\sum_{K,i}c^{Ki}_I\frac{\partial^2\varphi_K}{\partial x_i\partial r}+\hbox{first order}.
$$
The second inequality is then obtained by applying 
$\zeta_0\Lambda^{s-k-2}_{\rm tan}\frac{\partial^k}{\partial r^k}$ to the above
equation and taking $L_2$ norms.

Using elliptic regularization one sees that all partial derivatives of $\zeta_0\varphi$
are square-integrable and the theorem follows.\enddemo

\vglue9pt\centerline{\bf  Appendix}
\vglue12pt

Here we will recall briefly the notion of ideal finite type and the explicit expressions for subelliptic multipliers as
introduced in [K2] and [K3]; this material is also explained in [DK]. Given an $(n-1)$-tuple $\rho_1,\dots,\rho_{n-1}$ of
germs of $C^\infty$ functions at $x_o\in M$ we denote by  ${\cal  M}(\rho_1,\dots,\rho_{n-1})$ the \\$(n-1)\times2(n-1)$
matrix defined by
$$
{\cal  M}(\rho_1,\dots,\rho_{n-1})=\pmatrix{c_{11}&c_{12}&\ldots&c_{1,n-1}\cr
                 c_{21}&c_{22}&\ldots&c_{2,n-1}\cr
                 \vdots&\vdots&\ddots&\vdots\cr 
                 c_{n-1,1}&c_{n-1,2}&\ldots&c_{n-1,n-1}\cr
                 L_1\rho_1&L_2\rho_1&\ldots&L_{n-1}\rho_1\cr
                 L_1\rho_2&L_2\rho_2&\ldots&L_{n}\rho_2\cr
                 \vdots&\vdots&\ddots&\vdots\cr
                 L_1\rho_{n-1}&L_2\rho_{n-1}&\ldots&L_{n-1}\rho_{n-1}\cr}.
$$ 
\proclaim{Theorem} If $M$ is a pseudoconvex {\rm CR} manifold of dimension $2n-1$ and if $x_o\in M$ then the   
${\cal  I}_q^+(x_o)$ have the following properties\/{\rm :}\/

\vglue4pt
{\bf A.} ${\cal  I}_q^+(x_o)$ is an ideal{\rm .}

\vglue4pt
{\bf B.} The real radical of ${\cal  I}_q^+(x_o)${\rm ,} denoted by $^{\Bbb   R}\hspace{-0.7em}\sqrt{{\cal  I}_q^+(x_o)}${\rm ,} is
contained in
${\cal  I}_q^+(x_o)$. The ideal $^{\Bbb   R}\hspace{-0.7em}\sqrt{{\cal  I}_q^+(x_o)}$ consists of all germs $g$ such that
there exist 
$f\in{\cal  I}_q^+(x_o)$ and $m\in\Bbb   Z$ with $|g|^m\le |f|$. 

\vglue4pt
{\bf C.} ${\rm Det}^{n-q-1}{\cal  M}(\rho_1,\dots,\rho_{n-1})$ is the ideal generated by the 
$$(n-q-1)\times(n-q-1)$$  subdeterminants of ${\cal  M}(\rho_1,\dots,\rho_{n-1})$. If $\rho_i\in{\cal  I}_q^+(x_o)$ 
then 
$${\rm Det}^{n-q-1}{\cal  M}(\rho_1,\dots,\rho_{n-1})\subset{\cal  I}_q^+(x_o).$$    
\endproclaim

\numbereddemo{Definition}
The ideals ${\cal  I}_{q,k}^+(x_o)$ are defined by induction on $k$ as follows:
$${\cal  I}_{q,1}^+(x_o)=\ ^{\Bbb R}\hspace{-0.7em}\sqrt{{\rm Det}^{n-q-1}{\cal  M}(0,\dots,0)},$$
$${\cal  I}_{q,k+1}^+(x_o)=\ ^{\Bbb R}\hspace{-0.7em}\sqrt{\Bigl({\cal  I}_{q,k}^+(x_o),\  {\rm Det}^{n-q-1}{\cal  M}(\rho_1,\dots,\rho_{n-1}})\  {\rm with}\ \rho_i\in
{\cal  I}_{q,k}^+(x_o)\Bigr).$$ 
\enddemo

Here $\Bigl({\cal A}, {\cal  B}, {\cal  C},\cdots\Bigr)$ denotes the ideal generated by 
$\{{\cal A}\cup {\cal  B}\cup {\cal  C}\cup\cdots\}$. The ideals ${\cal  I}_{q,k}^-(x_o)$ are defined by setting
${\cal  I}_{q,k}^-(x_o)={\cal  I}_{n-q-1,k}^+(x_o)$. We set 
${\cal  I}_{q,k}(x_o)={\cal  I}_{q,k}^+(x_o)\cap{\cal  I}_{q,k}^-(x_o)$.

\demo{{R}emark} It then follows that ${\cal  I}_{q,k}^+(x_o)\subset{\cal  I}_q^+(x_o)$, 
${\cal  I}_{q,k}^-(x_o)\subset{\cal  I}_q^-(x_o)$, and ${\cal  I}_{q,k}(x_o)\subset{\cal  I}_q(x_o)$. Furthermore, ${\cal 
I}_{q,k}^+(x_o)\!\subset\!{\cal  I}_{q+1,k}^+(x_o)$ and  ${\cal  I}_{q,k}^-(x_o)\!\supset\!{\cal  I}_{q+1,k}^-(x_o)$, so that
${\cal  I}_{q,k}(x_o)={\cal  I}_{m,k}^+(x_o),$ where  $m=\min\{q,\ n-q-1\}$.

\numbereddemo{Definition}If $\Omega$ is a pseudoconvex domain in a hermitian manifold $X$ with a smooth boundary $M$
then 
$x_o\in M$ is of {\it finite ideal $q$\/{\rm -}\/type} (for $\bar\partial$) if $1\in{\cal  I}_{q,k}^+(x_o)$ for some $k$. If $M$
is a pseudoconvex CR manifold and $x_o\in M$ then $x_o\in M$ is of {\it finite ideal $q$\/{\rm -}\/type} (for $\bar\partial_b$) if
$1\in{\cal  I}_{q,k}(x_o)$ for some $k$.
\enddemo

Thus, both for domains and CR manifolds finite ideal $q$-type implies that subellipticity holds for $(p,q)$-forms. The
question is whether this condition is also necessary. For domains in two-dimensional manifolds, necessity was proved by
Greiner (see [G]). In this case the ideal finite type condition can be expressed in terms of commutators of vector fields.
The proof is easily generalized to the case when the matrix $(c_{IJ})$ is diagonalizable in a neighborhood of $x_o$. It
also can be easily generalized to CR manifolds for which both matrices $(c_{IJ})$ and $(c_{I'J'})$ can be diagonalized in a
neighborhood of $x_o$. In case the defining function for the boundary $M$ is real analytic the results of Diederich and
Forn\ae ss (see [DF]) were used in [K2] to prove that finite ideal type is equivalent to subellipticity. In the general case
Catlin proved that subellipticity is equivalent to finite D'Angelo type (see [C] and [D]). Thus the problem is to prove
that finite D'Angelo type is equivalent to finite ideal type.

\AuthorRefNames [Ch2]

\end{document}